\documentclass[final,3p,times, 10pt]{elsarticle}

\usepackage{lineno}     
\modulolinenumbers[5]   

\usepackage[utf8]{inputenc}
\usepackage[T1]{fontenc}
\usepackage{microtype} 

\usepackage{amsmath, amssymb, amsbsy, amsfonts}
\usepackage{mathtools}
\usepackage[ruled,vlined]{algorithm2e} 

\usepackage{graphicx}
\usepackage{xcolor}
\usepackage{tikz}
\usetikzlibrary{positioning, arrows.meta, calc, shapes.geometric, backgrounds, fit}
\usepackage{pgfplots}
\usepackage{booktabs}
\pgfplotsset{compat=1.18} 
\usepackage{subcaption}
\usepackage{pgfplots}
\usepgfplotslibrary{fillbetween}
\usepackage{bbm}

\usepackage{booktabs} 
\usepackage{siunitx}
\biboptions{sort&compress}

\usepackage[colorlinks=true, linkcolor=blue, citecolor=blue, urlcolor=blue]{hyperref}
\usepackage[nameinlink]{cleveref} 

\DeclareMathOperator*{\argmin}{arg\,min}

\newtheorem{problem}{Problem}

\journal{Journal}

\begin{document}

\begin{frontmatter}


    \title{Constraint-driven Optimization and Parametrization of Industrial NURBS Geometries via Neural Deformation Field}

    \author[1]{Federico Tamburlin}
    \author[1]{Giovanni Canali}
    \author[1]{Giuseppe Alessio D'Inverno}    
    \author[1,3]{Nicola Demo}    
    \author[2]{Andrea Mola}
    \author[1,3]{Gianluigi Rozza}
        
    \affiliation[1]{organization={SISSA, International School for Advanced Studies},
            city={Trieste},
            postcode={34136}, 
            country={Italy}}
    
    \affiliation[2]{organization={IMT School for Advanced Studies Lucca},
            city={Lucca},
            postcode={55100}, 
            country={Italy}}          
            
    \affiliation[3]{organization={FAST Computing Srl},
            city={Trieste},
            postcode={34136}, 
            country={Italy}}

    \begin{abstract}
    This work presents a differentiable framework for the parametrization and shape optimization of industrial CAD geometries represented by multi-patch NURBS surfaces.
    The method enables the deformation of complex CAD models through a physics-informed geometric parametrization, allowing direct morphing driven by physical constraints without the need to prescribe a predefined deformation strategy.
    A neural displacement field, implemented as a multi-layer perceptron acting on the NURBS control points, provides a compact parametrization of the admissible design space while preserving patch connectivity.\\
Global geometric quantities relevant to hydrostatic design, including displaced volume, wetted surface area and buoyancy centroid, are formulated as differentiable integral operators evaluated on the parametric domain. These quantities are computed through Gauss–Legendre quadrature combined with analytical B-spline derivatives for surface metric evaluation, allowing gradient propagation to the deformation parameters while limiting the computational overhead of automatic differentiation.\\
The proposed framework operates directly on CAD representations without intermediate mesh generation. Numerical experiments on a modified KVLCC2 hull demonstrate the capability of the method to satisfy competing hydrostatic constraints while producing smooth CAD-compatible geometries and showing stable convergence across multiple random initializations.
    \end{abstract}

    \begin{highlights}
\item Differentiable CAD-native pipeline using neural parametrizations to optimize NURBS geometry via integral constraints.
\item End-to-end gradient propagation from loss to control points — no adjoint equations needed.
\item CAD-compatible output (IGES) preserving smoothness and topology for industrial downstream processes.
\end{highlights}


\end{frontmatter}



\section{Introduction}
\label{sec:intro}

\subsection{Motivation}
Computer-Aided Design (CAD) models are the native representation of industrial geometries in Product Lifecycle Management (PLM)~\cite{srinivasan2011integration,sudarsan2005product} and Computer-Aided Manufacturing (CAM)~\cite{bi2020computer}. Free-form shapes are commonly described through Boundary Representation (B-Rep) models based on Non-Uniform Rational B-Splines (NURBS)~\cite{versprille1975computer}, which provide an exact parametrization of curves and surfaces through control points, weights, knot vectors, and polynomial degrees. Preserving the geometric and topological validity of this representation during deformation is essential for downstream analysis and manufacturing.

Traditional shape optimisation workflows~\cite{IvagnesDemoRozza2023IJNME,DemoTezzeleMolaRozza2021JMSE,villa2021effective} require the definition of a reduced set of design parameters capable of producing meaningful geometric variations while preserving the validity of the underlying shape representation. Classical deformation-based parametrizations, such as Free-Form Deformation (FFD)~\cite{sederberg1986free} and Radial Basis Function (RBF) morphing, are widely used in these pipelines because they can modify complex geometries without rebuilding the CAD model from scratch. They can act on different representations, including CAD control entities, surface discretizations, point clouds, and volumetric meshes, and have been successfully applied in industrial and naval shape optimization through geometric, mesh, and adaptive grid deformation strategies~\cite{wei2025hull,villa2021effective,li2022application,wang2024cfd}. However, their effectiveness strongly depends on the definition of the deformation structure, such as the FFD lattice, the RBF control points, the support radius, or the propagation strategy. Coarse deformation models may fail to capture localized features, whereas finer parametrizations increase the number of design variables and may introduce spurious high-frequency modes~\cite{MolaDemoTezzeleRozzaCH17}. Moreover, since these deformation maps are usually external to the native CAD parametrization, preserving geometric regularity, derivative information, and global CAD-level constraints is not automatic, especially when the quantities of interest are integral properties rather than pointwise displacements or mesh-based measures.

In naval hydrodynamics, several works have addressed hull-form optimization through advanced deformation strategies acting either on the geometry or directly on the computational mesh. \cite{wei2025hull} proposed a three-dimensional deformation strategy for hull-form optimization, while~\cite{villa2021effective} introduced an effective mesh deformation approach for hull shape design optimization based on CAD-consistent geometric modifications and volumetric mesh propagation. Similarly, a combined mesh deformation with adaptive refinement in hull-form design was introduced in~\cite{li2022application}, and~\cite{wang2024cfd} presented a CFD-based optimization framework in calm water using adaptive grid deformation. These contributions demonstrate the practical relevance of deformation-based parametrizations in naval applications.

Alongside classical mesh and geometry deformation techniques, data-driven and generative approaches have recently emerged as promising tools to construct compact design spaces for complex industrial shapes. In~\cite{Lee2024} authors argued that traditional CAD parameters may limit the exploration of large shape variations and topological changes, whereas generative networks can provide low-dimensional manifolds for design optimization. In the naval field, ~\cite{khan2023shiphullgan} proposed a deep convolutional generative model for parametric hull design, showing how generative networks can learn compact representations of admissible ship hull shapes.~\cite{bagazinski2023shipgen} introduced a diffusion-based model for parametric ship hull generation under multiple objectives and constraints, whereas more generally,~\cite{padula2024generative} focused on generative models for the deformation of industrial shapes subject to linear geometric constraints, with emphasis on model-order and parameter-space reductions. These approaches highlight a growing interest in replacing high-dimensional, manually defined shape parameters with reduced latent or generative representations. However, implicit and generative descriptions may reduce direct control over standard CAD entities and may require additional mechanisms to guarantee strict geometric validity.

The advent of Deep Learning and Physics-Informed Neural Networks has also promoted the use of Automatic Differentiation as a mechanism for differential equation solution~\cite{raissi2019physics}. This makes it possible to compute sensitivities without deriving problem-specific adjoint codes~\cite{Sun2023}, which are traditionally required in conventional gradient-based shape optimization workflows~\cite{jalui2025developing,mohebbi2020exact}. Differentiable programming further allows geometric kernels to be embedded directly into gradient-based loops.~\cite{Moola2024} proposed \textit{THB-Diff}, exploiting Automatic Differentiation on GPUs for Truncated Hierarchical B-splines fitting and refinement, mainly targeting computational efficiency in spline evaluation. Current PINN-based shape optimization frameworks~\cite{Tillmann2023,Zhang2024} often rely on coordinate projection strategies or problem-specific parametrizations, differing from the direct optimization of standard CAD control variables.

In this context, the present work follows a complementary direction. Rather than relying exclusively on mesh-dependent deformation pipelines~\cite{villa2021effective,li2022application,wang2024cfd}, stochastic or generative design manifolds~\cite{khan2023shiphullgan,bagazinski2023shipgen,padula2024generative}, or soft-constrained PINN formulations~\cite{Tillmann2023,Zhang2024}, it defines an analytical and differentiable map from NURBS control points to global geometric metrics. Automatic Differentiation is combined with numerical quadrature to enforce integral geometric properties directly at the CAD-parametrization level, rather than purely minimizing fitting errors or evaluating constraints on a generated mesh. This enables gradient-based optimization while preserving explicit control over industrial shape deformation and bypassing explicit mesh generation during the constraint-evaluation stage.

\subsection{Main novelties and contributions}
\begin{figure*}[htbp]
    \centering
    \begin{tikzpicture}[
        scale=0.75, transform shape, 
        font=\rmfamily\small, 
        >=Stealth,
        base_block/.style={
            rectangle, rounded corners=4pt, draw=black, thick, align=center,
            minimum width=2.5cm
        },
        blue_block/.style={
            base_block, fill=cyan!10, minimum height=2.4cm
        },
        grey_block/.style={
            base_block, fill=gray!20, minimum height=2.4cm, align=left
        },
        red_block/.style={
            base_block, fill=red!60, minimum height=2.4cm, 
            font=\rmfamily\bfseries\large
        },
        mlp_block/.style={
            rectangle, rounded corners=6pt, draw=black, thick, fill=blue!30,
            minimum width=2.5cm, minimum height=0.9cm, font=\rmfamily\bfseries
        },
        line/.style={draw, thick, ->},
        dashed_line/.style={draw, thick, dashed, ->}
    ]

        \node[red_block] (spline) {SPLINE\\model};

        \node[mlp_block, above=1.5cm of spline] (mlp) {MLP};

        \node[blue_block, left=1.5cm of spline] (domain) {
            \textbf{PARAMETRIC}\\
            \textbf{DOMAIN}\\[0.3cm]
            $(u,v) \in [0,1]^2$
        };

        \node[blue_block] (cad) at (domain |- mlp) {
            \textbf{\large CAD}\\[0.8cm] 
        };
        \begin{scope}[shift={($(cad.center)+(0,-0.2)$)}]
            \draw[thick, fill=blue!5] (0,0.15) ellipse (0.7 and 0.2);
            \draw[thick, fill=white, fill opacity=0.6] (-0.7,0.15) arc (180:360:0.7 and 0.55);
            \draw[thin, gray] (-0.7,0.15) .. controls (-0.3,-0.5) and (0.3,-0.5) .. (0.7,0.15);
            \draw[thin, gray] (0,-0.4) .. controls (0,-0.1) .. (0,0.35);
        \end{scope}

        \node[blue_block, right=1.5cm of spline] (physical) {
            \textbf{PHYSICAL}\\
            \textbf{DOMAIN} (x,y,z)\\[0.8cm]
        };
        
        \begin{scope}[shift={($(physical.center)+(0,-0.3)$)}, scale=0.55]
            \coordinate (A) at (-0.8, -0.6);
            \coordinate (B) at (0.5, -0.9);
            \coordinate (C) at (1.1, -0.2);
            \coordinate (D) at (-0.2, 0.2);
            \coordinate (E) at (-0.9, 0.7);
            \coordinate (F) at (0.7, 0.9);
            \draw[gray!60, thin] (E)--(D) (D)--(B);
            \draw[black!80, thin] (A)--(B)--(C)--(F)--(E)--(A);
            \draw[black!80, thin] (A)--(D)--(F);
            \draw[black!80, thin] (B)--(C);
            \foreach \p in {A,B,C,D,E,F} \fill[black] (\p) circle (3pt);
        \end{scope}

        \node[grey_block, right=1.5cm of physical, align=left] (loss) {
            \textbf{\large LOSS}\\[0.2cm]
            $L = \lambda_V L_{\text{VOL}} + \lambda_C L_{\text{CENTR}} +$\\
            $\quad \lambda_D L_{\text{DRAFT}} + \lambda_A  L_{\text{WSA}} +$\\
            $\quad \lambda_R L_{\text{REG}}$
        };

        \draw[line] (cad.south east) -- node[above right] {pyiges} (spline.135); 
        
        \draw[line] (mlp.south) -- (spline.north);
        
        \draw[dashed_line] (domain.east) -- node[midway, above=0.2cm, font=\scriptsize, align=center] {SAMPLING\\$(u_i, v_i) \; w_i$} (spline.west);
        
        \draw[line] (spline.east) -- (physical.west);
        \draw[line] (physical.east) -- (loss.west);

        \draw[thick, ->] (loss.north) 
            |- node[pos=0.75, above] {UPDATE} (mlp.east); 

\end{tikzpicture}  
    \caption{Schematic of the end-to-end differentiable optimization pipeline. The system maps original geometrical control points through a Neural Displacement Field (MLP) to generate the deformed surface. Geometric properties (Volume, Area, Centroid) are computed via analytical B-Spline derivatives and numerical quadrature, driving the loss function composed of integral constraints, geometric barrier and regularization terms.}
    \label{fig:arch_1}    
\end{figure*}
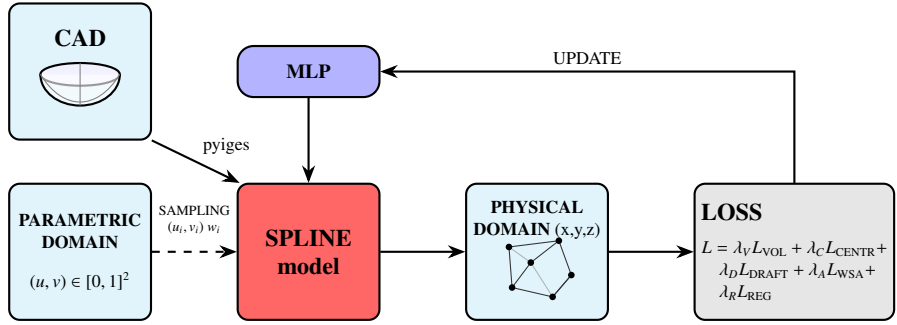

This work presents a shape design and optimization pipeline for parametric surfaces (\Cref{fig:arch_1}), implemented within the PINA framework \cite{coscia2023physics}. 
The proposed approach focuses on a simplified physical setting, where the quantities of interest are global geometric and hydrostatic measures, such as volume, surface area, and centroid-related properties, evaluated through analytical B-spline derivatives and numerical quadrature. While full PDE-constrained shape optimization remains the standard approach for many engineering applications, its integration within a differentiable CAD-native framework introduces significant complexity and is out of the scope of the present work, representing a natural direction for future extensions. Instead, this paper introduces a constraint-driven deformation paradigm, where the geometry is optimized with respect to physically meaningful integral targets, avoiding explicit field equation solves while retaining a clear physical interpretation. In this context, the framework can be interpreted within the broader class of PINNs, where physical consistency is enforced through differentiable constraints embedded in the loss function, and not with data. The deformation is modeled via a Neural Displacement Field, implemented as a Multi-Layer Perceptron (MLP), which maps the original CAD control points to their deformed configuration, enabling a compact and flexible parametrization of the design space.

The main contributions of this paper are summarized as follows:
\begin{itemize}
    \item \textbf{Neural Parametrization of Geometric Deformation:} An approach employing a Multi-layer Perceptron (MLP) to parametrize the displacement field of NURBS control points. This provides a flexible and compact representation of the design space compared to point-by-point optimization and establishes a global deformation law that maintains interface continuity \cite{hoyer2019neural, chandrasekhar2021tounn}.
    \item \textbf{End-to-End Differentiable Mapping:} By integrating the B-Spline mapping directly into the computational graph, the framework enables exact propagation of gradients from the objective function to the control points. This eliminates the need for manually derived adjoint equations.
    \item \textbf{Enforcement of Integral Geometric Constraints:} The formulation of differentiable loss terms to satisfy global integral constraints, such as target volume. This ensures that the optimized geometry adheres to functional global requirements, implemented within the PINA framework.
    \item \textbf{CAD Compatibility:} The pipeline operates directly on standard NURBS patch parameters, maintaining the topological validity and smoothness required for industrial downstream processes. The output is a standard CAD file (e.g., IGES).
\end{itemize}

\subsection{Outline of the paper}
The remainder of this paper is organized as follows.  \Cref{sec:preliminary} introduces the main theoretical concepts and notation. \Cref{sec:methodology} details the implementation within the PINA framework and the neural parametrization of deformation. \Cref{sec:results} presents numerical experiments validating the proposed approach on a multi-patch ship geometry. Finally, \Cref{sec:conclusions} summarizes the findings and outlines future possible research directions.


\section{Preliminary notation and concepts}
\label{sec:preliminary}

This section outlines the mathematical formulation of the geometric representation and the numerical integration scheme adopted within the differentiable optimization pipeline.

\subsection{B-Splines Surface Parametrization}
Let $\hat{\Omega} = [0,1]^2 \subset \mathbb{R}^2$ denote the parametric reference domain. A physical surface $\Omega \subset \mathbb{R}^3$ is considered, defined by the mapping $\mathbf{S}: \hat{\Omega} \to \Omega$. In the context of standard CAD formats, the geometry is represented as a Non-Uniform Rational B-Spline (NURBS) surface, which is a subset of the B-Spline employed in this work:
\begin{equation}
    \mathbf{S}(u,v) = \sum_{i=0}^{n} \sum_{j=0}^{m} R_{i,j}^{p,q}(u,v) \,\mathbf{C}_{i,j}, \quad (u,v) \in \hat{\Omega},
\end{equation}
where $\mathbf{C}_{i,j} \in \mathbb{R}^3$ are the control points forming the control net. The rational basis functions $R_{i,j}^{p,q}$ are defined as:
\begin{equation}
    R_{i,j}^{p,q}(u,v) = \frac{N_{i,p}(u) N_{j,q}(v) w_{i,j}}{\sum_{k=0}^{n} \sum_{l=0}^{m} N_{k,p}(u) N_{l,q}(v) w_{k,l}},
\end{equation}
where $N_{i,p}(u)$ and $N_{j,q}(v)$ are the univariate B-spline basis functions of degree $p$ and $q$, defined over the knot vectors $\Xi_u = \{\xi_0, \dots, \xi_{n+p+1}\}$ and $\Xi_v = \{\eta_0, \dots, \eta_{m+q+1}\}$, respectively, and $w_{i,j}$ are the associated weights \cite{piegl2012nurbs}. 

In the proposed framework, the topological parameters, specifically the knot vectors, polynomial degrees, and weights, are fixed to the initial geometry values. The spatial coordinates of the control points $\mathbf{C} = \{\mathbf{C}_{i,j}\}$ constitute the set of design variables subject to optimization.

\subsection{Geometry and Metric Distortion}
Evaluation of geometric quantities on the physical manifold $\Omega$ requires the computation of the metric distortion induced by the mapping $\mathbf{S}$. The tangent basis vectors are given by the partial derivatives $\partial_u \mathbf{S}$ and $\partial_v \mathbf{S}$. The local surface area element $d\Omega$ is related to the parametric area element $d\hat{\Omega}$ via the Jacobian determinant $J$, defined as the magnitude of the normal vector field:
\begin{equation}
    J(u,v) = \|\mathbf{n}(u,v)\|_2, \quad \text{with} \quad \mathbf{n}(u,v) = \frac{\partial \mathbf{S}}{\partial u} \times \frac{\partial \mathbf{S}}{\partial v}.
\end{equation}
This Jacobian term ensures that gradients computed in the parametric space correctly reflect the sensitivity of physical objective functionals in the Euclidean space. The computation of the tangent vectors relies on the exact analytical derivatives of the univariate B-spline basis functions. For a basis function $N_{i,p}(u)$ of degree $p$ defined on a knot vector $U$, the derivative is calculated recursively as:
\begin{equation}
    \frac{d}{du} N_{i,p}(u) = \frac{p}{u_{i+p}-u_i} N_{i,p-1}(u) - \frac{p}{u_{i+p+1}-u_{i+1}} N_{i+1,p-1}(u).
\end{equation}

\subsection{Integral Constraints and Numerical Quadrature}
Global geometric properties, such as volume and centroids, are evaluated via numerical integration over the parametric domain. A Gauss-Legendre \cite{abramowitz1964handbook} quadrature rule is employed with $N_q$ nodes $\{(\hat{u}_k, \hat{v}_k)\}_{k=1}^{N_q}$ and corresponding weights $\{\omega_k\}_{k=1}^{N_q}$. The approximation of a generic scalar field integral $I = \int_{\Omega} f(\mathbf{x}) \, d\Omega$ is given by:
\begin{equation}
    I \approx \sum_{k=1}^{N_q} f(\mathbf{S}(\hat{u}_k, \hat{v}_k)) \, J(\hat{u}_k, \hat{v}_k) \, \omega_k.
\end{equation}

\subsubsection{Volumetric and Surface Constraints}
For hydrodynamic applications, the geometry is often modelled as an open surface bounded by a reference plane at $z = z_{\text{ref}}$. The enclosed volume $V$ is computed by integrating the vertical elevation weighted by the projected area on the $xy$-plane. Leveraging the divergence theorem, this is formulated as:
\begin{equation}
    V = \int_{\Omega} (z - z_{\text{ref}}) \, n_z \, d\Omega \approx \sum_{k=1}^{N_q} (S_z(\hat{u}_k, \hat{v}_k) - z_{\text{ref}}) \, n_z(\hat{u}_k, \hat{v}_k) \, J(\hat{u}_k, \hat{v}_k) \, \omega_k,
\end{equation}
where $n_z$ is the $z$-component of the unit normal. In the numerical implementation, the term $n_z J$ is computed efficiently as the $z$-component of the un-normalized cross product $\partial_u \mathbf{S} \times \partial_v \mathbf{S}$, avoiding explicit normalization steps.

The Center of Buoyancy (COB) is evaluated directly on the displaced volume by computing the static moments of the submerged geometry. Leveraging the divergence theorem on the bounding surface, the spatial coordinates are formulated as:

\begin{equation}
\mathbf{x}_{COB} = \frac{\int_{\Omega} \mathbf{X} (z - z_{ref}) n_z \, d\Omega}{\int_{\Omega} (z - z_{ref}) n_z \, d\Omega}, \quad \mathbf{X} = \begin{bmatrix} x \\ y \\ \frac{z + z_{ref}}{2} \end{bmatrix}.
\label{eq:cob_vector}
\end{equation}

\noindent
where the term $\frac{z + z_{\text{ref}}}{2}$ represents the vertical lever arm of the elementary volume column.

Finally, the Wetted Surface Area (WSA) corresponds to the surface area of the manifold:
\begin{equation}
    A = \int_{\Omega} \mathbbm{1}_{z \le z_{\text{ref}}} \, d\Omega \approx \sum_{k=1}^{N_q} \mathbbm{1}_{z(\hat{u}_k, \hat{v}_k) \le z_{\text{ref}}} \, J(\hat{u}_k, \hat{v}_k) \, \omega_k.
\end{equation}
The evaluation of these quantities and their differentiability is central to the proposed optimization pipeline.


\section{Methodological Framework}
\label{sec:methodology}

Consider a computational domain $\Omega \subset \mathbb{R}^3$ defined by a set of NURBS patches, in the physical space. The shape optimization problem is formulated as the search for a deformation field that satisfies global integral and geometrical constraints while minimizing a regularity functional. The proposed framework relies on a global neural parametrization of the displacement field, integrated into a differentiable computing environment.

\subsection{Neural Deformation Operator}
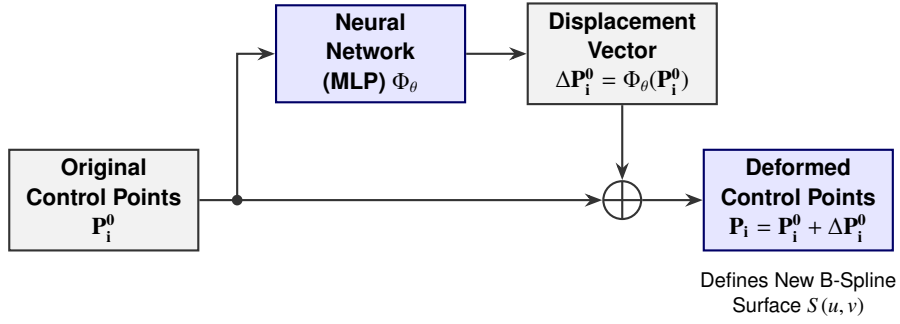
\begin{figure*}
\centering
    \begin{tikzpicture}[
    node distance=0.6cm and 1.0cm, 
    >=Stealth,
    basebox/.style={
        rectangle, 
        draw=black!80, 
        thick, 
        minimum width=2.5cm, 
        minimum height=1.1cm, 
        align=center, 
        font=\sffamily\bfseries\small, 
        fill=gray!10
    },
    bluebox/.style={basebox, draw=blue!40!black, fill=blue!10},
    sum/.style={
        circle, 
        draw=black!80, 
        thick, 
        minimum size=0.5cm, 
        inner sep=0pt,
        path picture={
            \draw[black!80, thick] (path picture bounding box.south) -- (path picture bounding box.north);
            \draw[black!80, thick] (path picture bounding box.west) -- (path picture bounding box.east);
        }
    },
    line/.style={
        draw=black!80, 
        thick, 
        ->
    }
]

    \node[basebox] (orig) {Original\\Control Points\\$\mathbf{P_i^0}$};    
    \node[bluebox, above right=0.6cm and 1.0cm of orig] (mlp) {Neural\\Network\\(MLP) $\Phi_{\theta}$};    
    \node[basebox, right=0.8cm of mlp] (disp) {Displacement\\Vector\\$\Delta \mathbf{P^0_i} = \Phi_{\theta}(\mathbf{P^0_i})$};    
    \node[sum] (plus) at (orig.center -| disp.center) {};    
    \node[bluebox, right=0.8cm of plus] (deformed) {Deformed\\Control Points\\$\mathbf{P_i} = \mathbf{P_i^0} + \Delta \mathbf{P^0_i}$};    
    \node[below=0.15cm of deformed, font=\sffamily\footnotesize, align=center] {Defines New B-Spline\\Surface $S(u,v)$};

    \path (orig.east) -- (orig.east -| mlp.west) coordinate[midway] (splitpoint);
    \draw[line] (orig.east) -- (plus.west);    
    \draw[line] (splitpoint) |- (mlp.west);    
    \draw[line] (mlp.east) -- (disp.west);    
    \draw[line] (disp.south) -- (plus.north);    
    \draw[line] (plus.east) -- (deformed.west);
    \fill[black!80] (splitpoint) circle (2pt);
\end{tikzpicture}
    \caption{Neural Deformation architecture: the displacement is computed by the MLP and added to the original control points.}
    \label{fig:mlp_arch}
\end{figure*}

Let $\mathcal{C}^0 = \{ \mathbf{P}^0_i \}_{i \in \mathcal{I}}$ be the set of control points defining the reference geometry $\Omega_0$. Let us define a continuous deformation map $\Phi_{\boldsymbol{\theta}}: \mathbb{R}^3 \to \mathbb{R}^3$, parametrized by a Multi-Layer Perceptron (MLP) with trainable weights $\boldsymbol{\theta}$ (\Cref{fig:mlp_arch}). 
To ensure a balance between expressivity and geometric regularity, the architecture comprises 5 hidden layers with hyperbolic tangent activation, followed by a linear projection to the physical space $\mathbb{R}^3$. This depth was selected based on preliminary benchmarks: deeper networks ($>8$ layers) tended to introduce optimization instability and unwanted high-frequency artifacts, reducing the implicit regularization effect. Conversely, shallower networks ($<3$ layers) lacked the spectral capacity to capture localized shape variations, such as the deformation required at the bulbous bow.
The set of deformed control points $\mathcal{C}(\boldsymbol{\theta}) = \{ \mathbf{P}_i(\boldsymbol{\theta}) \}_{i \in \mathcal{I}}$ is defined by applying the displacement field to the reference set $\mathcal{C}^0$:

\begin{equation}
    \mathbf{P}_i(\boldsymbol{\theta}) = \mathbf{P}_i^0 + \Phi_{\boldsymbol{\theta}}(\mathbf{P}_i^0), \quad \forall i \in \mathcal{I}.
    \label{eq:neural_map}
\end{equation}

Since $\Phi_{\boldsymbol{\theta}}$ is a deterministic smooth function ($C^\infty$), topological connectivity at patch interfaces is satisfied by construction. For any shared boundary $\Gamma_{AB} = \Omega_A \cap \Omega_B$, coincident control points receive identical displacement vectors, preserving the water-tightness of the assembly.
To enforce symmetry with respect to the plane $y=0$, the operator output is explicitly constrained. The hard constraint is imposed:

\begin{equation}
    \Phi_{\boldsymbol{\theta}, y}(\mathbf{x}) = 0, \quad \forall \mathbf{x} \text{ s.t. } |y| < \epsilon.
\end{equation}

Analogously, a hard constraint is imposed on the vertical displacement of the upper boundary. Letting $z_{top}$ denote the maximum vertical coordinate of the reference hull, the operator is constrained such that:
\begin{equation}
    \Phi_{\boldsymbol{\theta},z}(\mathbf{x})= 0, \quad \forall \mathbf{x} \text{ s.t. } |z - z_{top}| < \epsilon.
\end{equation}

The network $\Phi_{\boldsymbol{\theta}}$ is initialized with a distribution ensuring zero initial displacements, allowing the optimization to start from the exact reference configuration $\Omega_0$.

\subsection{Integral Constraints and Geometric Measures}
The shape optimization objectives are formulated as global integral quantities computed over the parametric domain. The PINA framework is extended to evaluate global geometric integrals via numerical quadrature.
A hybrid differentiation scheme is implemented: the neural network parameters are optimized via standard automatic differentiation, while the surface geometric properties (normals and Jacobian determinants) are computed using exact analytical B-Spline basis derivatives. This avoids the explicit automatic differentiation of the B-spline kernel for tangent vector computation, reducing the computational graph overhead by approximately 50\% \footnote{Benchmarks on the Jacobian evaluation per epoch show a reduction from 0.365\unit{s} (Automatic Differentiation) to 0.193\unit{s}  (Analytical PINA derivative), confirming a computational saving of $\approx47$\%}.

Let $\mathcal{Q}(\Omega)$ be a generic quantity of interest defined as a rational integral:
\begin{equation}
    \mathcal{Q}(\Omega) = \frac{\int_{\hat{\Omega}} N(\mathbf{u}(\boldsymbol{\xi})) J(\boldsymbol{\xi}) d\boldsymbol{\xi}}{\int_{\hat{\Omega}} D(\mathbf{u}(\boldsymbol{\xi})) J(\boldsymbol{\xi}) d\boldsymbol{\xi}} \approx \frac{\sum_{k} N_k J_k w_k}{\sum_{k} D_k J_k w_k},
    \label{eq:rational_integral}
\end{equation}
where $\mathbf{u}$ is the displacement field, and $J_k$ is the geometric measure derived from the analytical B-Spline derivatives. Specifically, $J = \|\mathbf{n}\|$ is used for surface metrics and $J = \mathbf{n} \cdot \mathbf{e}_z$ for volumetric quantities via the Divergence Theorem. The specific constraints are listed in \Cref{tab:integral_constraints}.

\begin{table}[htbp]
    \centering
    \caption{Definition of Integral Constraints formulation. Operators $N$ and $D$ refer to \Cref{eq:rational_integral}. For entries marked with $\emptyset$, the denominator is set to unity (scalar integral). Note: $\mathbf{X}$ is the vector $[x, y, \frac{z+z_{ref}}{2}]$.}
    \label{tab:integral_constraints}
    \small
    \begin{tabular}{lcccc}
        \toprule
        Constraint & $N(\mathbf{x})$ & $D(\mathbf{x})$ & Measure $J$ & Notes \\
        \midrule
        Displaced Volume & $z - Z_{ref}$ & $\emptyset$ & $\mathbf{n} \cdot \mathbf{e}_z$ & Projection on $xy$-plane \\
       Wetted Surface & $\mathbbm{1}_{z \le z_{ref}}$ & $\emptyset$ & $\|\mathbf{n}\|$ & Submerged area ($z \le z_{ref}$) \\
        Center of Buoyancy & $\mathbf{X} (z - z_{ref})$ & $z - z_{ref}$ & $\mathbf{n} \cdot \mathbf{e}_z$ & Volume moments \\
        \bottomrule
    \end{tabular}
\end{table}

\subsection{Optimization Problem Statement}
The optimization is formulated as an unconstrained minimization problem of a composite scalar loss function $\mathcal{L}(\boldsymbol{\theta})$. 

\begin{problem}[Shape Optimization]
\label{prob:shape_opt}
Find the optimal network parameters $\boldsymbol{\theta}^* \in \mathbb{R}^d$ such that:
\begin{equation}
    \boldsymbol{\theta}^* = \argmin_{\boldsymbol{\theta}} \left( \mathcal{L}_{int}(\boldsymbol{\theta}) + \mathcal{L}_{barrier}(\boldsymbol{\theta}) + \mathcal{L}_{reg}(\boldsymbol{\theta}) \right),
\end{equation}
subject to the hard constraints on symmetry and vertical displacement bounds.
\end{problem}



\noindent The \textit{Integral Constraints} $\mathcal{L}_{int}$, computed as the deviation from target geometric properties, (Volume $V$, Area $A$ and Center of Buoyancy $x_{COB}$) are minimized via a Sum of Relative Squared Errors:

\begin{equation}
\mathcal{L}_{int} = \lambda_{V} \left( \frac{V(\boldsymbol{\theta}) - V_{targ}}{V_{targ}} \right)^2+ \lambda_{A} \left( \frac{A(\boldsymbol{\theta}) - A_{targ}}{A_{targ}} \right)^2+ \lambda_{C} \left( \frac{x_{COB}(\boldsymbol{\theta}) - x_{COB,targ}}{x_{COB,targ}} \right)^2.
\end{equation}

\noindent A specific component mask is applied to the COB condition to penalize exclusively the longitudinal shift ($x_{COB}$), allowing the transverse and vertical components ($y_{COB}$, $z_{COB}$) to naturally adapt to the new displaced mass distribution.

The \textit{Geometric Barrier} $\mathcal{L}_{barrier}$ constraint, a unilateral constraint on the draft (minimum vertical coordinate $Z_{min}$), is enforced via a penalty method. The violation set is defined, based on the Rectified Linear Unit (ReLU) operator:

\begin{equation}
    \mathcal{L}_{barrier} = \lambda_{D} \sum_{k=1}^{N_q} \left[ 
    \max(0, Z_{min} - z(\mathbf{x}_k(\boldsymbol{\theta})))\right]^2.
\end{equation}

\noindent A displacement-dependent energy functional is introduced, the \textit{Regularization Functional} $\mathcal{L}_{reg}$, to ensure surface quality and suppress high-frequency geometric artifacts. Unlike Euclidean regularization on parameters, this term acts on the physical displacement $\mathbf{u}(\boldsymbol{\xi})$ numerically integrated over the manifold:

\begin{equation}
    \mathcal{L}_{reg} = \lambda_{reg} \sum_{k=1}^{N_q} \left[ \exp\left( \gamma_{reg} \| \mathbf{u}(\boldsymbol{\xi}_k) \|^2 \right) - 1 \right] J(\boldsymbol{\xi}_k) w_k.
\end{equation}

\noindent The exponential kernel imposes an asymptotic penalty as displacements approach the soft limit defined by $\gamma_{reg}^{-1/2}$, implicitly constraining the search space to small deformations, under a certain tunable threshold.

\subsection{Optimization Scheme}

The continuous problem is discretized via Gaussian quadrature \cite{abramowitz1964handbook}. The discrete optimization procedure is detailed in \Cref{alg:opt_loop}. The iterative scheme is structured into four primary stages. First, the forward pass of the neural network computes the unconstrained displacement field $\Delta \tilde{\mathcal{C}}$, which is subsequently projected onto the feasible space via the operator $\mathcal{P}$ to strictly enforce symmetry and boundary conditions.

Second, to minimize computational overhead, a geometric caching strategy is implemented. The physical coordinates $\mathbf{x}_q$ and the metric vectors $\mathbf{m}_q$ are evaluated exclusively at the invariant quadrature nodes. This approach isolates the exact evaluation of the B-spline analytical derivatives from the automatic differentiation graph associated with the neural network parameters.

Third, the functional assembly aggregates the geometric residuals $\mathcal{R}$ and the regularization energy $\mathcal{E}_{reg}$ via numerical integration. Finally, the network weights $\boldsymbol{\theta}$ are updated. The abstract operator \text{Optimizer} relies on the exact gradients $\nabla_{\boldsymbol{\theta}} \mathcal{L}$ extracted via automatic differentiation, utilizing AdamW for initial non-linear exploration, followed by the L-BFGS algorithm for second-order local refinement.

\begin{algorithm}[htbp]
\SetAlgoLined
\DontPrintSemicolon
\caption{Optimization Loop}
\label{alg:opt_loop}

\KwIn{Reference grid $\mathcal{C}^0$, Quadrature rule $(\Xi, W)$.}
\KwOut{Optimized parameters $\boldsymbol{\theta}^*$.}

Initialize MLP weights $\boldsymbol{\theta}_0$ \;
Compute reference configuration $\mathbf{x}^0(\Xi)$\;
\BlankLine
\While{$t < T_{\max}$}{
    \tcp{Neural Deformation Field}
    $\Delta \tilde{\mathcal{C}} \gets \Phi_{\boldsymbol{\theta}_t}(\mathcal{C}^0)$\;
    $\Delta \mathcal{C} \gets \mathcal{P}(\Delta \tilde{\mathcal{C}})$\;
    $\mathcal{C} \gets \mathcal{C}^0 + \Delta \mathcal{C}$\;
    
    \BlankLine
    \tcp{Geometric Analysis}
    \For{$q=1, \dots, N_q$}{
         $\mathbf{x}_q \gets \mathbf{S}(\boldsymbol{\xi}_q; \mathcal{C})$\;
         $\mathbf{t}_{u,q} \gets \partial_u \mathbf{S}(\boldsymbol{\xi}_q; \mathcal{C})$\;
         $\mathbf{t}_{v,q} \gets \partial_v \mathbf{S}(\boldsymbol{\xi}_q; \mathcal{C})$\;
         $\mathbf{m}_q \gets \mathbf{t}_{u,q} \times \mathbf{t}_{v,q}$\;
    }
    
    \BlankLine
    \tcp{Functional Assembly}
    $\mathcal{R}_V \gets \left( \frac{\sum_q (z_q - z_{ref}) (\mathbf{m}_q \cdot \mathbf{e}_z) w_q - V_{target}}{V_{target}} \right)^2$\;
    $\mathcal{R}_A \gets \left( \frac{\sum_{z_q \le z_{ref}} \|\mathbf{m}_q\| w_q - A_{target}}{A_{target}} \right)^2$\;
    $\mathcal{R}_{C} \gets \left( \frac{\sum_q x_q (z_q - z_{ref}) (\mathbf{m}_q \cdot \mathbf{e}_z) w_q}{x_{COB,target} \sum_q (z_q - z_{ref}) (\mathbf{m}_q \cdot \mathbf{e}_z) w_q} - 1 \right)^2$\;
    $\mathcal{R}_{D} \gets \sum_q \left( \max(0, Z_{min} - z_q) \right)^2$\;
    $\mathcal{E}_{reg} \gets \sum_q \left( e^{\gamma \|\mathbf{x}_q - \mathbf{x}^0_q\|^2} - 1 \right) \|\mathbf{m}_q\| w_q$\;
    
    \BlankLine
    \tcp{Parameter Update}
    $\mathcal{L}(\boldsymbol{\theta}_t) \gets \lambda_V \mathcal{R}_V + \lambda_A \mathcal{R}_A + \lambda_C \mathcal{R}_C + \lambda_{D}\mathcal{R}_{D} + \lambda_{reg}\mathcal{E}_{reg}$\;
    $\boldsymbol{\theta}_{t+1} \gets \text{Optimizer}(\boldsymbol{\theta}_t, \nabla_{\boldsymbol{\theta}} \mathcal{L}(\boldsymbol{\theta}_t))$\;
}
\end{algorithm}


\section{Results}
\label{sec:results}

\subsection{Integral convergence analysis}
The accurate evaluation of global geometric constraints necessitates a numerical integration scheme. A Gauss-Legendre \cite{abramowitz1964handbook} quadrature rule is applied over the parametric domain $\hat{\Omega}=[0,1]^{2}$ to approximate surface and volumetric integrals.

\begin{figure}[htbp]
    \includegraphics[width=1\linewidth]{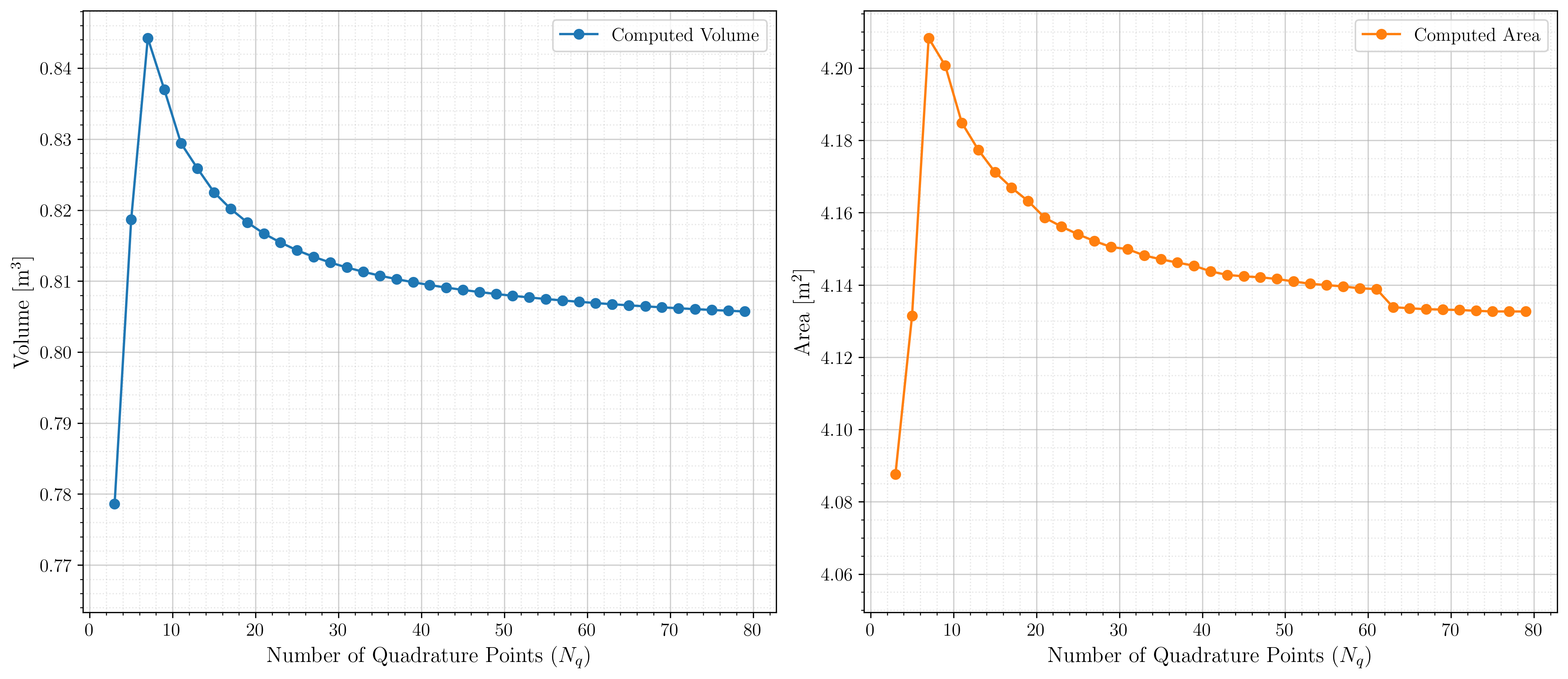}
    \caption{Convergence and stability analysis of geometric integrals via Gauss-Legendre quadrature. \\
    \textbf{Left panel:} Volumetric stability.
    \textbf{Right panel:} Convergence of the Wetted Surface Area (WSA).
    Both metrics motivate the selection of $N_q = 18$ per dimension as a sufficient trade-off resolution, ensuring that numerical integration noise remains negligible compared to the analytical gradients driving the neural deformation field.}
    \label{fig:convergence_analysis_2}
\end{figure}

As illustrated in \Cref{fig:convergence_analysis_2}, the volumetric computation exhibits convergence behaviour. Despite these initial error, the scheme achieves stabilisation.

For the optimization framework, a resolution of $N_{q} = 18$ points per dimension was established as a trade-off.


\subsection{Experimental Setup and Initial Benchmark}
The proposed framework was validated on a modified version of the \textit{KVLCC2} hull geometry\footnote{The input geometry was parsed using the \textit{pyiges} library \cite{cad-pinn-sources}.}  \cite{cad-pinn-sources}. 
While the framework is extensible to trimmed topologies, this experimental setup leverages the native tensor-product structure of untrimmed patches.
The optimization pipeline was executed on a consumer-grade mobile workstation\footnote{Hardware specifications: Intel Core i7 9th Gen H-series, 16 GB DDR4 RAM, NVIDIA GeForce GTX 1650 Mobile; OS: Ubuntu 24.04 LTS.}.

The optimization is performed using a fixed random seed to ensure reproducibility. The MLP weights are initialized using Kaiming initialization for hidden layers and zero-initialization for the output layer, to start from an identity mapping, with zero initial deformation of control points.\\
The numerical validation is conducted on a modified version of KVLCC2 \cite{cad-pinn-sources}
and all results are reported in model scale.
The primary objective of this numerical experiment was to assess the solver's capability to handle competing geometric constraints. The optimization problem was formulated to reduce the Wetted Surface Area (WSA) and the maximum Draft (minimum $z$ coordinate) while strictly conserving the volume and the longitudinal position of the Center of Buoyancy (COB).
The specific targets were defined relative to the initial geometric properties:
\begin{itemize}
    \item Volume Conservation: $V_{target} = 1.0 \cdot V_0$
    \item Wetted Area Minimization: $A_{target} = 0.97 \cdot A_0$ \quad (3\% reduction)
    \item Draft Constraint: $Z_{min}^{target} = 0.95 \cdot Z_{min}^{0}$ \quad (5\% reduction in depth)
    \item Centroid Preservation: ${x}_{COB}^{target} = {x}_{COB}^{0}$
\end{itemize}
This setup represents a preliminary design scenario where a vessel's efficiency is improved (lower friction drag via area reduction) and operational flexibility is increased (lower draft), without sacrificing capacity (volume).

\begin{table}[htbp]
    \centering
    \caption{Optimization parameters and configuration settings for the \textit{KVLCC2 modified} optimization.}
    \label{tab:hyperparameters}
    \begin{tabular}{lc}
        \toprule
        \textbf{Parameter} & \textbf{Value} \\
        \midrule
        Input File & \texttt{kvlcc2m.iges} \\
        Seed & 911\\
        Reference Plane ($Z_{ref}$) [\unit{m}] & 0.0 \\
        Quadrature Points ($N_q$) & 18 \\
        Optimization Epochs & 2000 \\
        \midrule
        \textbf{MLP Architecture} & \\
        Hidden Neurons & 64 \\
        Hidden Layers & 5 \\
        Learning Rate & $1 \cdot 10^{-3}$ \\
        \hline
        \textbf{Loss Weights} & \\
        $\lambda_{vol}$ (Volume) & 2.0 \\
        $\lambda_{area}$ (Area) & 1.0 \\
        $\lambda_{depth}$ (Draft Barrier) & 1.0 \\
        $\lambda_{cen}$ (Centroid) & 1.0 \\
        $\lambda_{reg}$ (Regularization) & $1 \cdot 10^{-4}$ \\
        $\gamma_{reg}$ (Exp. Decay) & 200.0 \\
        \bottomrule
    \end{tabular}
\end{table}

\subsubsection{Numerical Convergence}
The optimization was performed over 2000 epochs and
\Cref{tab:hyperparameters} details the optimization parameters utilized.
The processed CAD model represents a symmetric half-hull. To preserve the tensor-product structure required for exact numerical quadrature, the geometry utilizes the native untrimmed NURBS patches, intentionally omitting the minor deck trimming present in the official KVLCC2 benchmark, having a slightly different area and volume compared to the original hull.

\begin{table*}[htbp]
    \centering
    \caption{Optimization results for the modified KVLCC2 (seed 911). Volume and Area errors are relative [\%]. Centroid deviations are in \unit{meters}.}
    \label{tab:results}
    \sisetup{
        table-format=-1.5,
        table-number-alignment=center, 
        group-digits=false
    }
    \begin{tabular}{
        l
        S
        S
        S
        r
    }
    \toprule
    \textbf{Metric} & {\textbf{Initial}} & {\textbf{Target}} & {\textbf{Final}} & {\textbf{Error}} \\
    \midrule
    Volume [\unit{m^3}]       & 0.80165 & 0.80165 & 0.80104 & -0.076 \% \\
    Wetted Area [\unit{m^2}] & 4.12894 & 4.00507 & 4.01145 & +0.159 \% \\
    Max Draft ($z_{\min}$) [\unit{m}] & -0.35873 & -0.34079 & -0.34085 & -0.018 \% \\
    \midrule
    \multicolumn{5}{l}{\textit{Absolute Geometric Deviations}} \\
    \midrule
    Centroid $x$ [\unit{m}]  & -0.22528 & -0.22528 & -0.22527 & 0.00001 [\unit{m}] \\
    Centroid $y$ [\unit{m}]  &  0.22813 &  0.22813 &  0.23973 & 0.01160 [\unit{m}] \\
    Centroid $z$ [\unit{m}]  & -0.17117 & -0.17117 & -0.16060 & 0.01057 [\unit{m}] \\
    Centroid Shift ($L_2$) [\unit{m}] & {--} & {--} & {--} & 0.01569 [\unit{m}] \\
    \bottomrule
    \end{tabular}
\end{table*}

As presented in \Cref{tab:results}, the barrier constraint was satisfied within an absolute algorithmic tolerance ($1e-3$). The final minimum $z$ value reached $-0.34085$, deemed acceptable.
The solver achieved an area of approximately $4.011$, reaching the target of $4.005$ (approximately $3\%$ reduction). The final relative error against the target was $+0.159\%$, corresponding to a net physical reduction of $2.845\%$ from the initial state.
Volume conservation was achieved with high accuracy, resulting in a final deviation of $-0.076\%$. This indicates that the solver successfully balanced the volume loss induced by the draft reduction through deformation in other regions.

Regarding the Center of Buoyancy (COB), the longitudinal position ($x$) was preserved with negligible error ($0.00001$~m). However, a total 3D shift of $0.01569$~m was observed, primarily driven by vertical ($z$) and transverse ($y$) displacements. This shift is consistent with the geometric transformation required to maintain constant volume while reducing draft: as the bottom is raised, volume is redistributed laterally and vertically, naturally elevating the vertical center of buoyancy. The observed total shift in the total center of buoyancy does not compromise the hydrostatic equilibrium of the hull. The optimization of a symmetric semi-geometry strictly enforces transverse symmetry ($y=0$), preventing any induced heel. Furthermore, the conservation of the immersed volume ensures that the design waterline remains consistent, independent of the variation in the vertical center of buoyancy ($z_{COB}$).

\begin{figure}[h]
    \centering    \includegraphics[width=\linewidth]{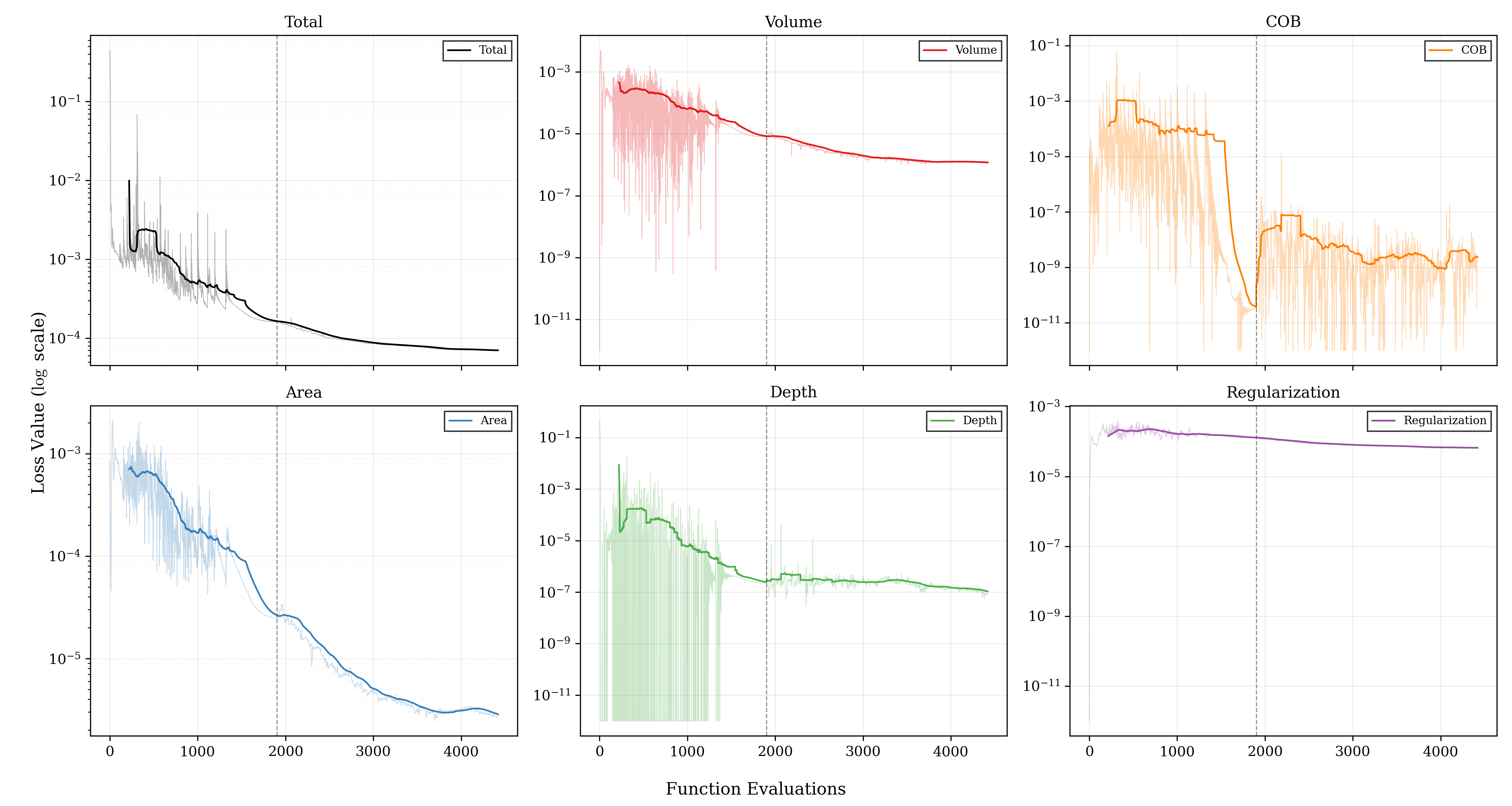}
    \caption{Optimization convergence analysis (Seed$=911$ and $\lambda_{reg}=1e^{-4}$). The plots illustrate the evolution of the global loss and specific geometric constraint losses against the number of function evaluations. Raw oscillation data is shown with high transparency, overlaid by a solid line representing a moving average to enhance trend visibility. The vertical dashed line indicates the switch from the stochastic exploration phase (AdamW) to the deterministic second-order refinement phase (L-BFGS). While the optimization was executed over 2000 epochs, the accumulation of more than 4000 function evaluations is a direct consequence of the L-BFGS algorithm, which inherently requires multiple evaluations per iteration for line search routines.}
    \label{fig:loss_911}
\end{figure}

\subsubsection{Geometric Deformation Analysis}
\label{sec:deform}

\begin{figure*}[htbp]
    \centering
    \begin{minipage}[t]{0.49\textwidth}
        \centering
        \includegraphics[width=\textwidth]{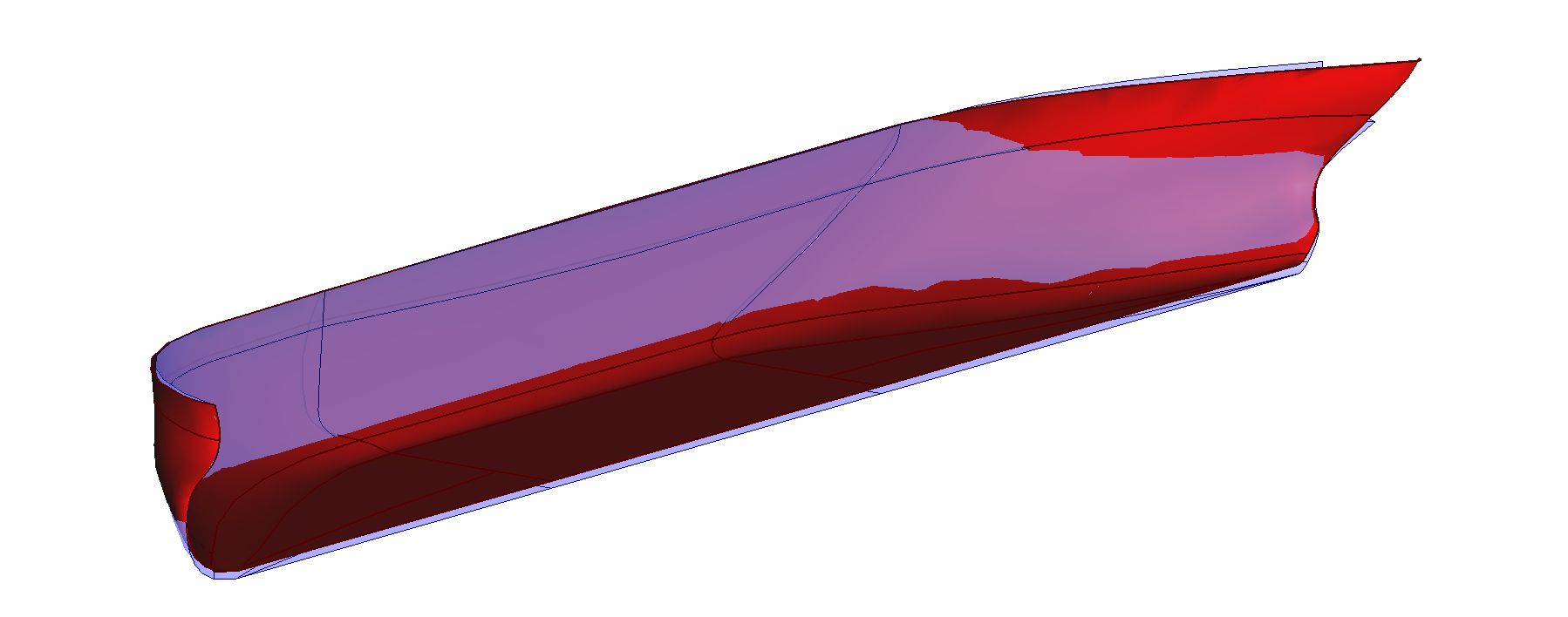}
    \end{minipage}
    \hfill
    \begin{minipage}[t]{0.49\textwidth}
        \centering
        \includegraphics[width=\textwidth]{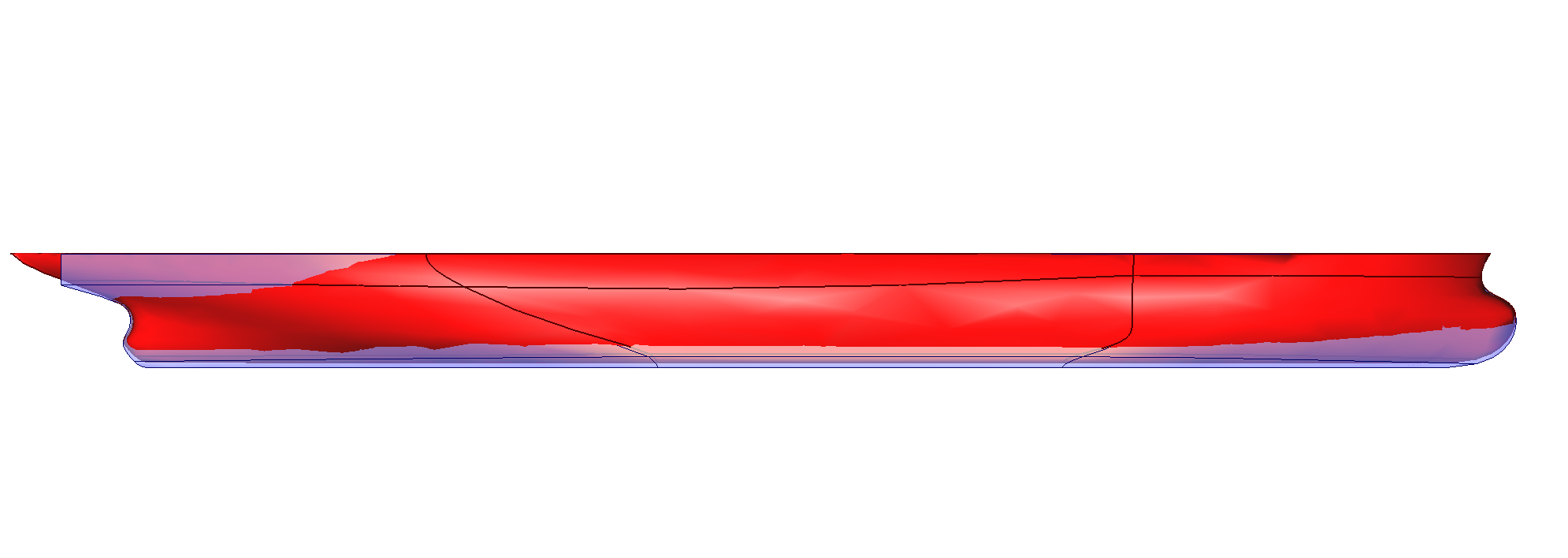}
    \end{minipage}
    \caption{
        Visualization of the optimized geometric deformation. Overlay of the initial hull (blue) and the optimized hull (red). Note the significant beam expansion at the midship section and the reshaping of the bulbous bow, generated by the network to compensate for the volume loss induced by the strict draft reduction constraint ($Z_{min}$). (Seed $= 911$ and $\lambda_{reg}=1e^{-4}$).
    }
    \label{fig:cad_comparison_3d_split}
\end{figure*}

\Cref{fig:cad_comparison_3d_split} shows the render comparison between the original surface and the deformed one.
The solver satisfied the draft constraint by effectively modifying the bottom nodes. To counteract the volume loss caused by the reduced draft, the neural network generates a transverse expansion the hull beam slightly in the midship section and reshape the bulbous bow, as suggested by the increase in the $y$-component of the centroid. This non-rigid deformation demonstrates the capability of the MLP-based parametrization to discover complex shape compensations.

\subsection{Benchmark against Free-Form Deformation (FFD)}
\label{sec:ffd_benchmark}

Free-Form Deformation (FFD) represents a standard methodology for CAD shape optimization and is utilized as the baseline to evaluate the proposed neural parametrization. The Spline architecture is still used, but the MLP deformation net is replaced by a FFD. This is implemented as a differentiable module computing the displacement field via Bernstein polynomials over a global spatial bounding box that encapsulates the NURBS control points \cite{sederberg1986free}. An anisotropic grid resolution of $24 \times 8 \times 8$ is utilized, introducing 4608 independent degrees of freedom. To ensure a rigorous comparative analysis, the objective function, geometric barrier, integral constraints, regularization penalties, and the dual stage optimization scheme remain identical to the MLP configuration.

Optimization metrics for the FFD baseline are detailed in \Cref{tab:results_ffd}. The FFD solver achieves numerical convergence, satisfying the minimum draft barrier with a relative error of -0.012\% and conserving the displaced volume to within -0.036\%. 

\begin{table}[htbp]
    \centering
    \caption{Optimization results for the modified KVLCC2 using Free Form Deformation (seed 911). Volume and Area errors are relative [\%].} 
    \label{tab:results_ffd}
    \sisetup{
        table-format=-1.5,
        table-number-alignment=center, 
        group-digits=false
    }
    \begin{tabular}{
        l
        S
        S
        S
        r
    }
    \toprule
    \textbf{Metric} & {\textbf{Initial}} & {\textbf{Target}} & {\textbf{Final}} & {\textbf{Error}} \\
    \midrule
    Volume [\unit{m^3}]       & 0.80165 & 0.80165 & 0.80136 & -0.036 \% \\
    Wetted Area [\unit{m^2}] & 4.12894 & 4.00507 & 4.00909 & +0.100 \% \\
    Max Draft ($z_{\min}$) [\unit{m}] & -0.35873 & -0.34079 & -0.34083 & -0.012 \% \\
    \midrule
    \multicolumn{5}{l}{\textit{Absolute Geometric Deviations}} \\
    \midrule
    Centroid $x$ [\unit{m}]  & -0.22528 & -0.22528 & -0.22527 & 0.00000 [\unit{m}] \\
    Centroid $y$ [\unit{m}]  &  0.22813 &  0.22813 &  0.23848 & +0.01035 [\unit{m}] \\
    Centroid $z$ [\unit{m}]  & -0.17117 & -0.17117 & -0.16190 & +0.00926 [\unit{m}] \\
    Centroid Shift ($L_2$) [\unit{m}] & {--} & {--} & {--} & 0.01389 [\unit{m}] \\
    \bottomrule
    \end{tabular}
\end{table}

Furthermore, the FFD methodology imposes an intrinsic resolution trade-off between geometric fidelity and computational cost: coarse FFD grids lack the local support necessary to capture and deform high-frequency geometric features, such as the bulbous bow, whereas dense grids exponentially increase the degrees of freedom. This unregularized behavior is incompatible with the objectives of this work, as the resulting geometric degradation compromises the topological validity and surface smoothness required by industrial standards. Unlike the non-linear MLP mapping, which adapts simultaneously to global structural modifications and local morphing requirements without grid dependencies, the extrinsic FFD grid is agnostic to the underlying topology. Consequently, it may more easily induce optimization stagnation when the required geometric transformations do not align with the principal axes of the control lattice.

\begin{figure*}[htbp]
    \centering
    \begin{subfigure}[b]{0.49\textwidth}
    \includegraphics[width=\textwidth]{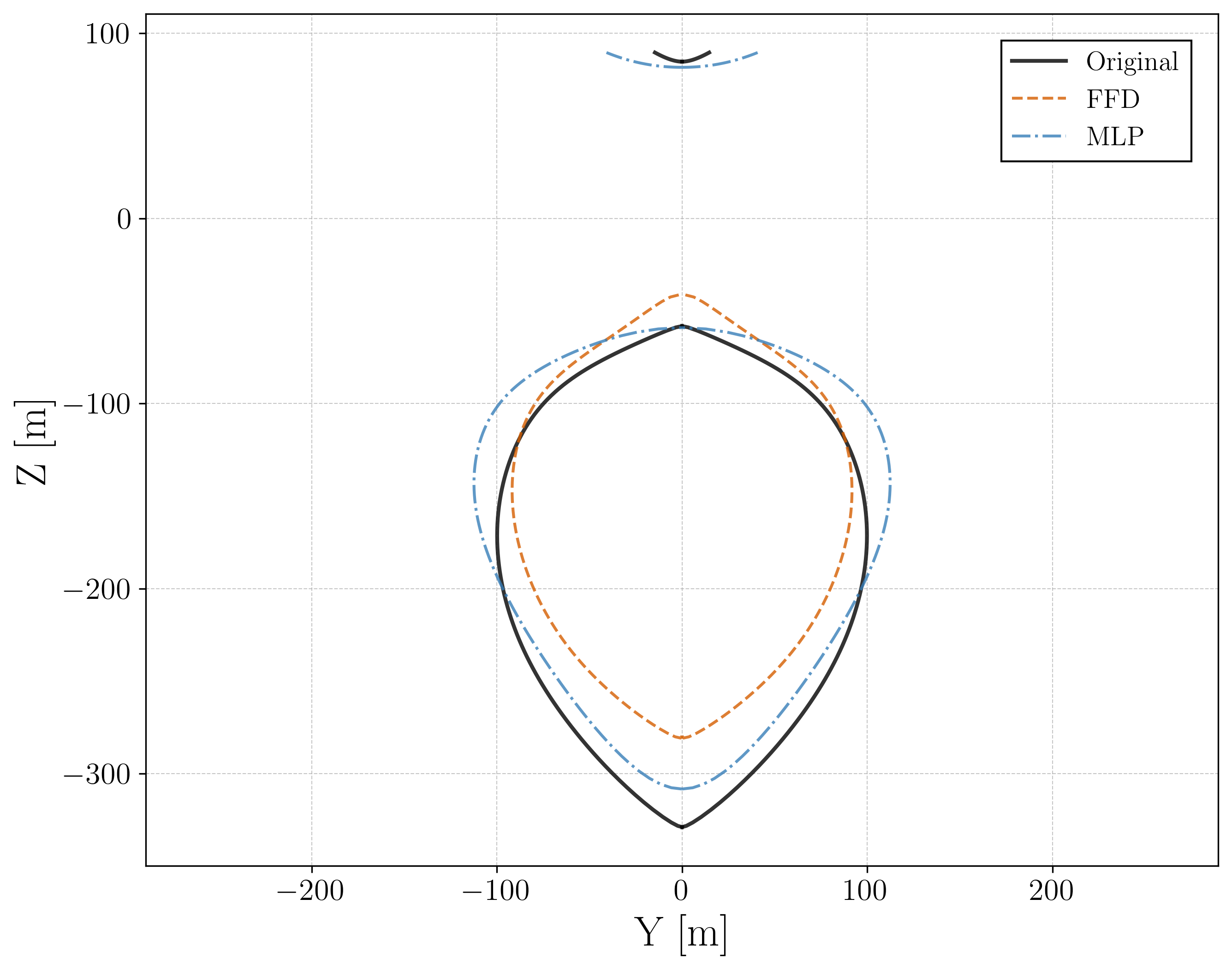}
        \caption{Bow Section (1/20)}
        \label{fig:bow_section_ffd}
    \end{subfigure}
    \hfill
    \begin{subfigure}[b]{0.49\textwidth}   
    \includegraphics[width=\textwidth]{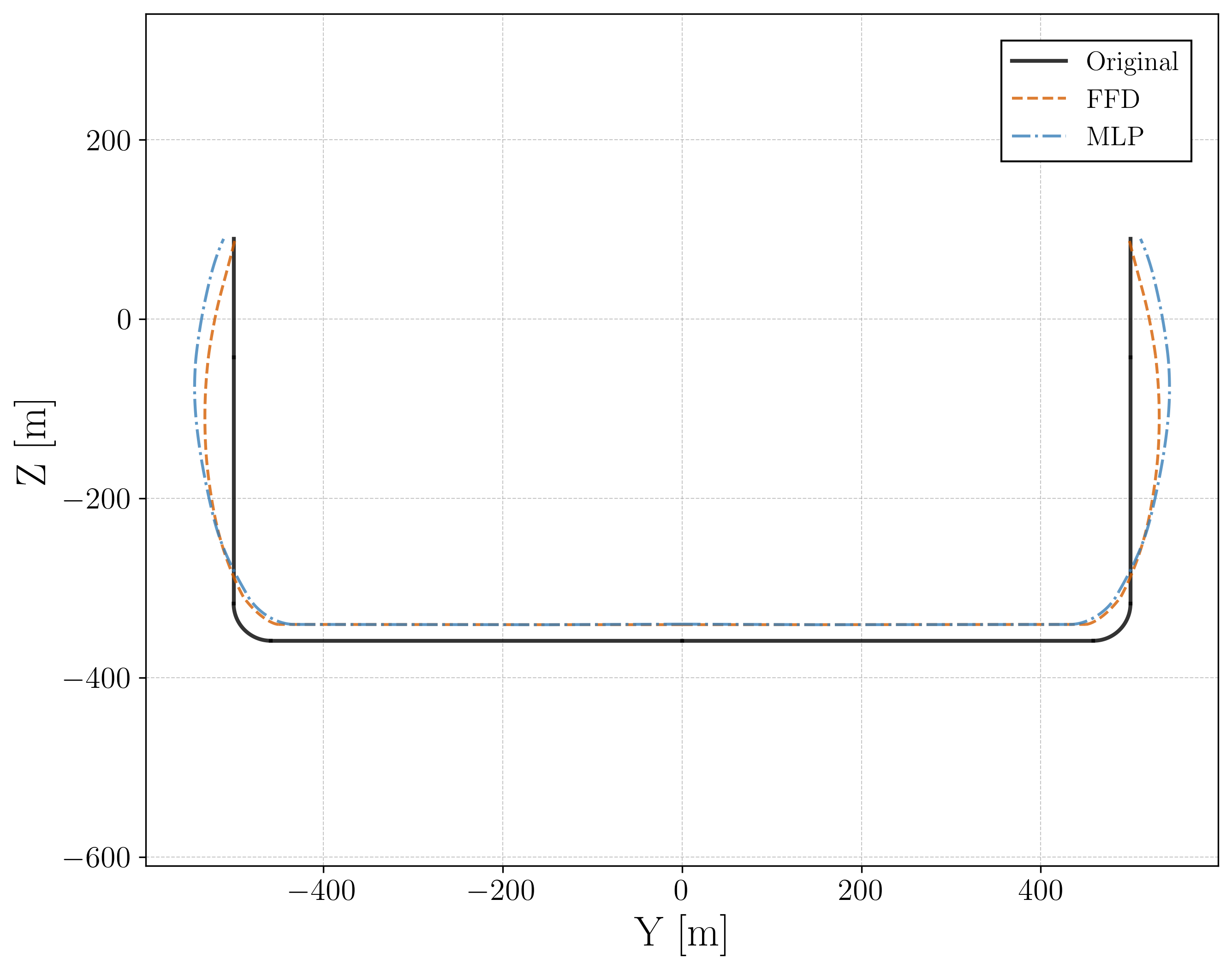}
        \caption{Midship Section (11/20)}        \label{fig:midship_section_ffd}
    \end{subfigure}
    \caption{Comparative cross-sectional analysis of the modified KVLCC2 hull geometry optimized via Free-Form Deformation (FFD) and the proposed neural parametrization (MLP). Same optimisation parameters.}
    \label{fig:ffd_mlp_sec}
\end{figure*}

 As observed in \Cref{fig:ffd_mlp_sec}, both methodologies produce comparable global deformation trends, in particular in the midship region. However, the MLP architecture exhibits superior flexibility in resolving complex shape compensations. 
 This capability is explicitly evident in the bow section (\Cref{fig:bow_section_ffd}), where the neural parametrization yields a more distinctive geometric profile compared to the conventional deformation induced by the FFD approach.



\subsection{Benchmark against Direct Control Points Optimization}

To quantify the regularization effect of the MLP Neural parametrization, a baseline experiment was conducted where the control points $\mathbf{C}_{i,j}$ were optimized directly as independent variables, bypassing the neural parameterization ($\Phi_{\boldsymbol{\theta}}$) from the pipeline. The optimization parameters are the same as above. While this approach allows for a theoretical zero-residual solution due to the maximized degrees of freedom, the resulting geometry exhibited severe topological degradation.

\begin{figure*}[htbp]
    \centering
    \begin{minipage}[t]{0.49\textwidth}
        \centering        \includegraphics[width=\textwidth]{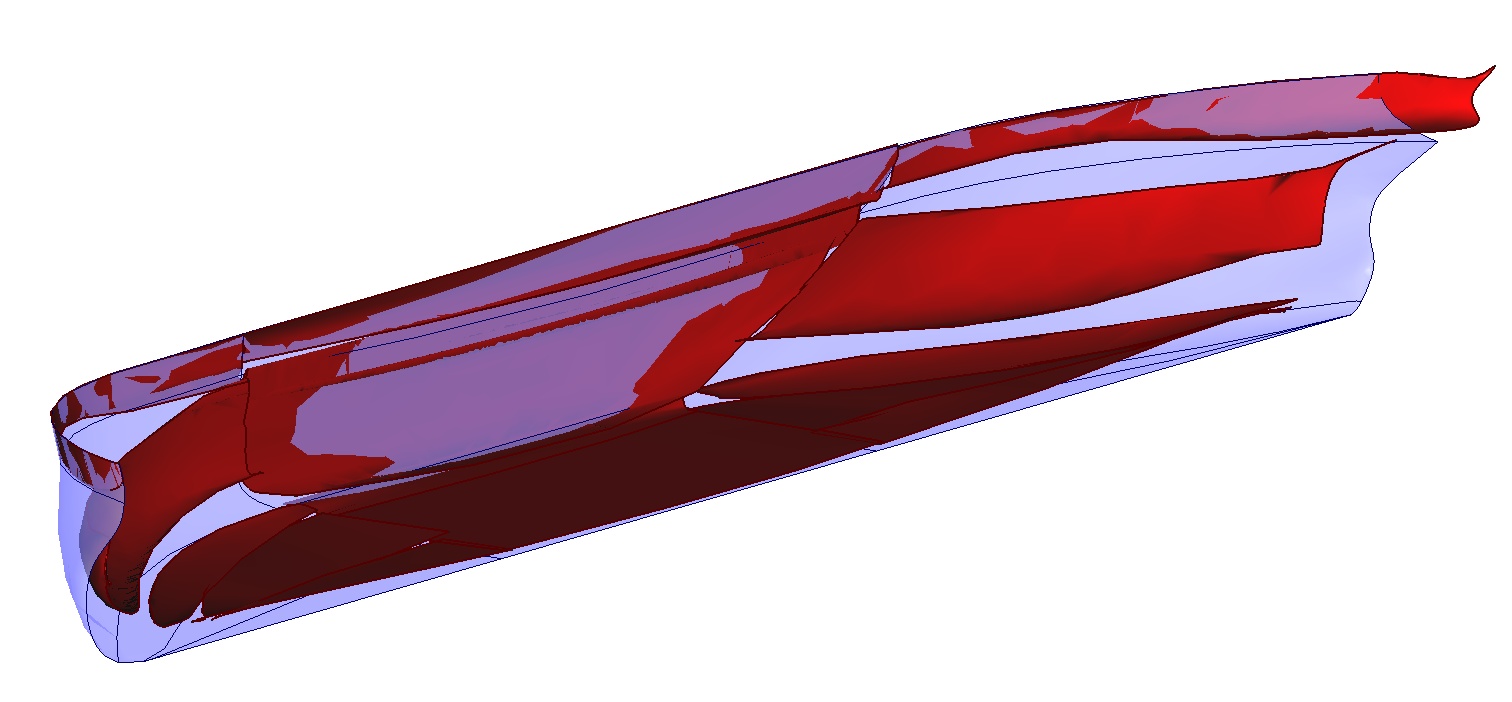}
    \end{minipage}
    \hfill
    \begin{minipage}[t]{0.49\textwidth}
        \centering        \includegraphics[width=\textwidth]{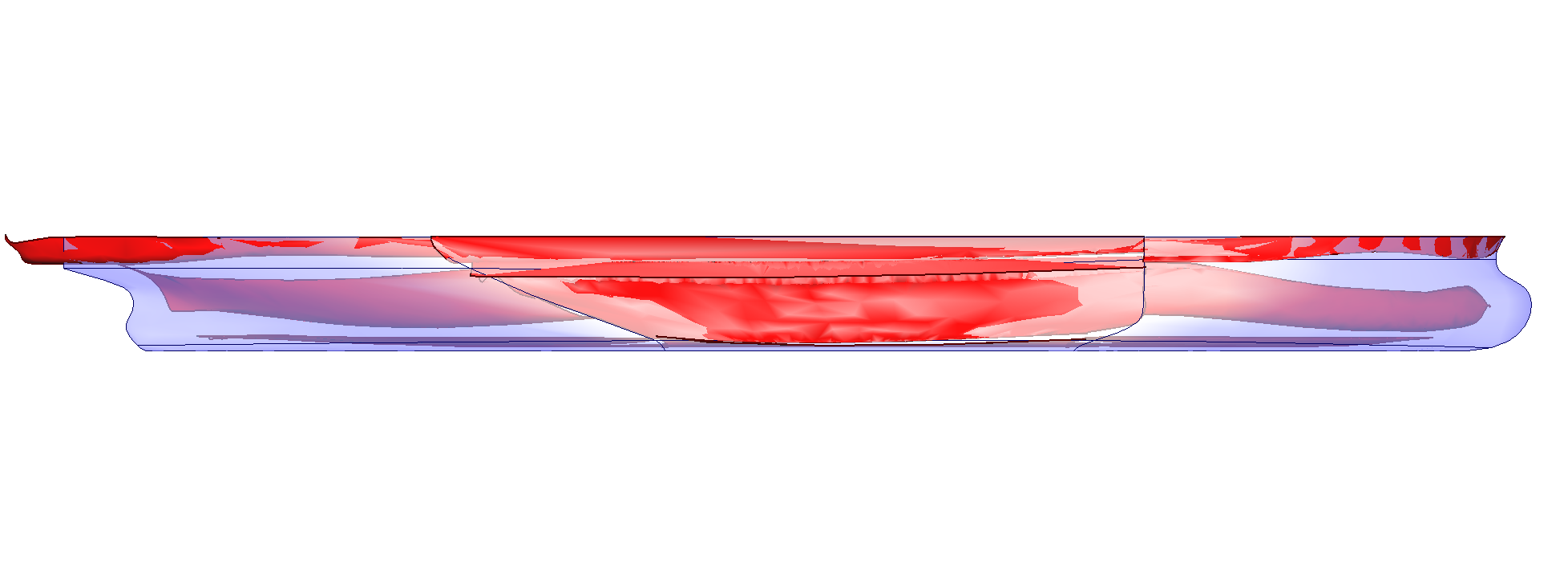}
    \end{minipage}
    \caption{
    Impact of the parametrization strategy on geometric quality. The direct optimization of control points lacks regularization, resulting in high-frequency artifacts and loss of $C^0$ continuity at patch interfaces (visible as surface gaps). The blue geometry is the original undeformed surface.
    }    
    \label{fig:direct_optimization}
\end{figure*}

\begin{figure*}[htbp]
    \centering
    \begin{subfigure}[b]{0.49\textwidth}
        \centering    \includegraphics[width=\textwidth]{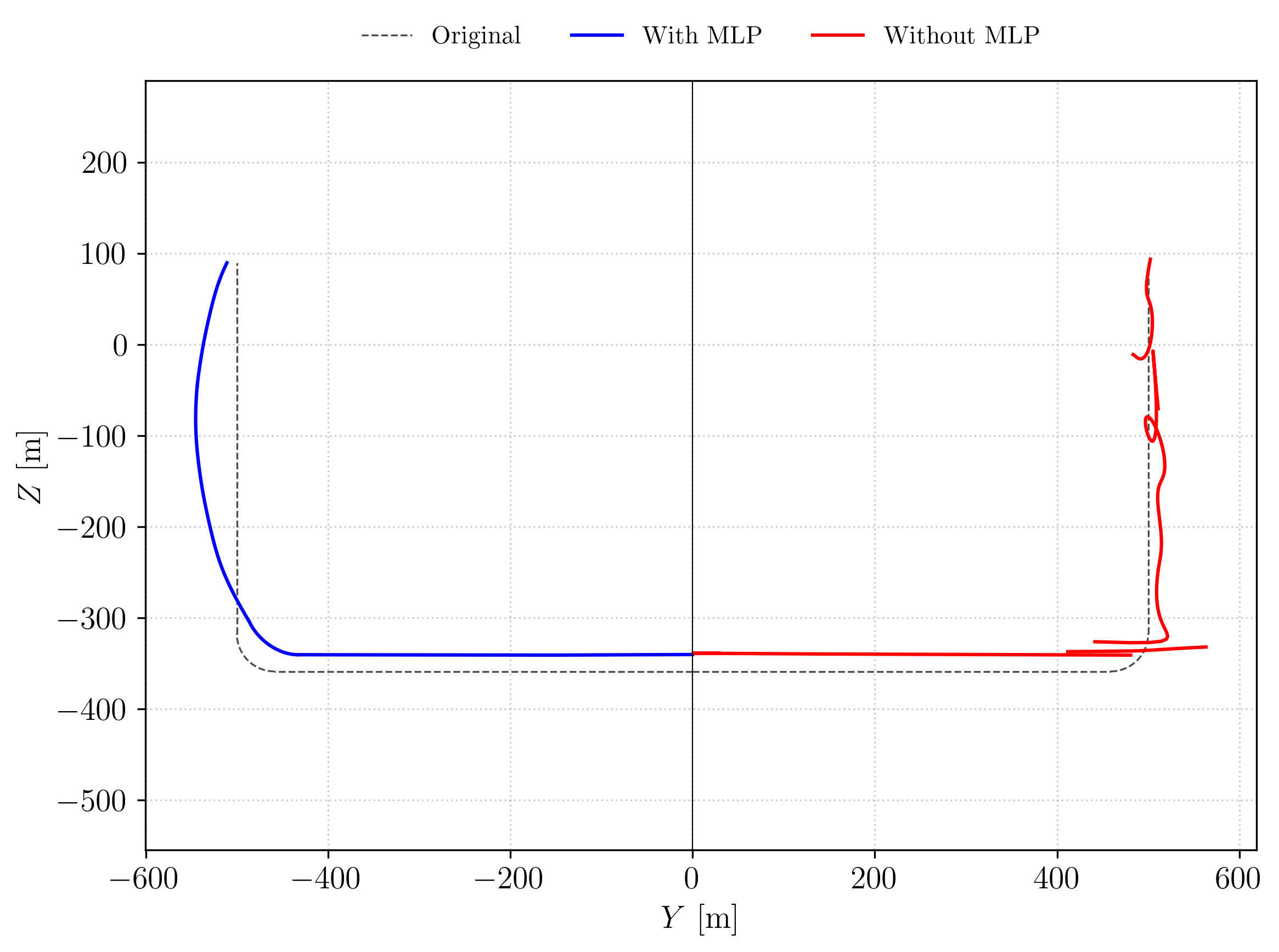}
        \caption{Midship Station (10/20)}
        \label{fig:midship_station}
    \end{subfigure}
    \hfill
    \begin{subfigure}[b]{0.49\textwidth}
        \centering          \includegraphics[width=\textwidth]{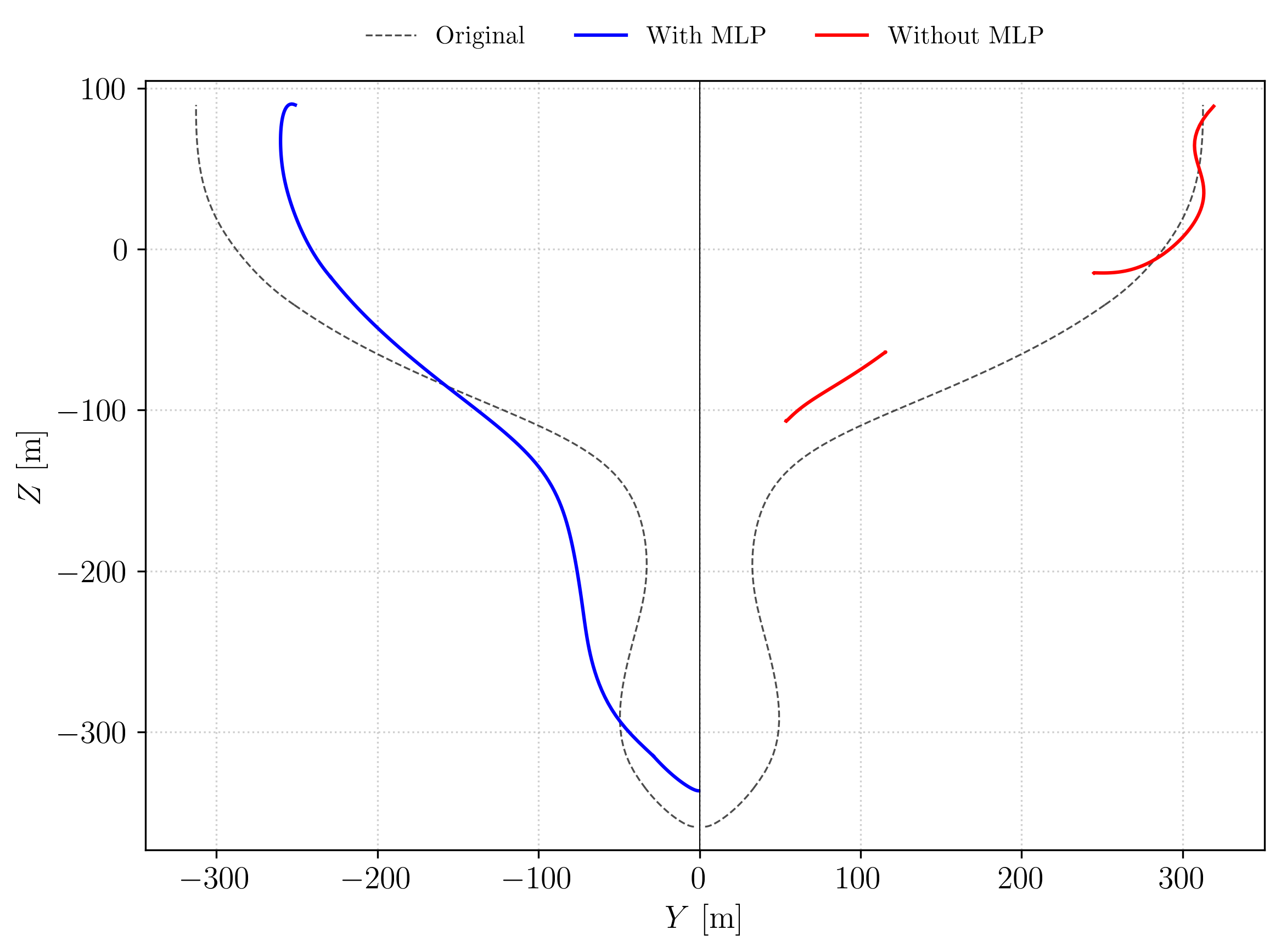}
        \caption{Aft Section (18/20)}
        \label{fig:aft_section}
    \end{subfigure}
    \caption{Comparative cross-sectional analysis at Station 10 and Station 18 over 20 total sections. The negative $Y$ domain ($Y<0$) visualizes the proposed MLP-driven deformation (blue), demonstrating strict adherence to topological constraints and the preservation of global hull continuity. In contrast, the positive $Y$ domain ($Y>0$) depicts the unregularized direct control point optimization (red), characterized by the breakdown of $C^0$ continuity and the macroscopic detachment of adjacent NURBS patches relative to the reference hull geometry (dashed grey).}
    \label{fig:mlp_sec}
\end{figure*}

As illustrated in \Cref{fig:direct_optimization} and \Cref{fig:mlp_sec}, the direct optimization lacks the regularisation effect of the neural network. Consequently, the control points drifted incoherently, resulting in high-frequency geometric noise and, critically, the violation of $C^0$ continuity at the patch interfaces. In the absence of the global deformation field provided by $\mathcal{N}_\theta$, maintaining watertight integrity would require the explicit imposition of computationally expensive equality constraints for every shared boundary node. 

While direct optimization treats each control point as an independent degree of freedom, making it prone to high-frequency artifacts, the MLP parametrization inherently favors low-frequency, smooth solutions due to the spectral bias of neural networks~\cite{rahaman2019spectral}. This effectively correlates the movement of spatially adjacent control points without requiring explicit smoothing penalties, consistent with findings in structural optimization~\cite{hoyer2019neural}.

Furthermore, the neural parametrization offers an advantage in terms of dimensionality. In a direct optimization setting, the number of degrees of freedom scales linearly with the number of control points, which can reach thousands for high-fidelity industrial models ($N_{dof} \approx 3 \times N_{cp}$). Conversely, the MLP parameter space is fixed by the network architecture ($N_{weights}$), independent of the geometric resolution, allowing for better control over memory consumption.


\subsection{Robustness Analysis and Variance}

To assess the dependence of the optimization outcome on the stochastic initialization of the neural network weights, the framework was executed using 20 distinct random seeds, using the same other optimisation parameters of \Cref{tab:hyperparameters}. \Cref{tab:different_seeds} summarizes the final geometric properties for each run.

\begin{table*}[htbp]
    \centering
    \caption{Hull optimization results. Statistical performance aggregated over 20 random seeds compared to \textit{Initial }and \textit{Target} values. Reported values represent Mean $\pm$ Standard Deviation.}
    \label{tab:different_seeds}
    
    \sisetup{
        table-number-alignment=center,
        group-digits=false,
        separate-uncertainty=true, 
        multi-part-units=single    
    }
    
    \begin{tabular}{
        l
        S[table-format=-1.5]      
        S[table-format=-1.5]      
        S[table-format=-1.5(6)]   
    }
    \toprule
    \textbf{Metric} 
        & {\textbf{Initial}}
        & {\textbf{Target}}
        & {\textbf{Statistical Analysis} (20 Seeds)} \\
    \midrule
    
    Final Volume [\unit{m^3}]        
        & 0.80165 & 0.80165 & 0.80086 +- 0.00064 \\
        
    Volume Error [\%]                
        & {--}    & 0.0     & -0.099 +- 0.079 \\
    \midrule
    
    Wetted Area [\unit{m^2}]         
        & 4.12894 & 4.00507 & 4.01630 +- 0.01268 \\
        
    Area Error [\%]                  
        & {--}    & 0.0     & +0.280 +- 0.317 \\
        
    Physical Area Delta [\%]         
        & {--}    & -3.0    & -2.728 +- 0.307 \\
    \midrule
    
    Max Draft ($z_{\min}$) [\unit{m}] 
        & -0.35873 & -0.34079 & -0.34090 +- 0.00003 \\
        
    Draft Dev. (from Limit) [\%]      
        & {--}     & 0.0      & -0.032 +- 0.010 \\
    \midrule
    
    \multicolumn{4}{l}{\textit{Absolute Geometric Deviations}} \\
    \midrule
    
    Centroid Shift ($L_2$) [\unit{m}] 
        & {--}    & 0.0     & 0.01566 +- 0.00112 \\
        
    Centroid Shift ($x$) [\unit{m}]    
        & {--}    & 0.0     & 0.00001 +- 0.00001 \\
        
    \bottomrule
    \end{tabular}
    
    \medskip
    \footnotesize
    \raggedright
    \noindent\textit{Metric Definitions.}
    Volume and wetted area errors are computed as relative deviations from targets. 
    The draft deviation is evaluated relative to the lower bound ($z_{\min} \ge z_{\mathrm{limit}}$); 
    negative values indicate a slight violation of the strict bound, tolerated by a small numerical relaxation of the depth barrier formulation.
\end{table*}

\begin{figure*}[htbp]
    \centering
    \begin{subfigure}[b]{0.49\textwidth}
        \centering
        \includegraphics[width=\textwidth]{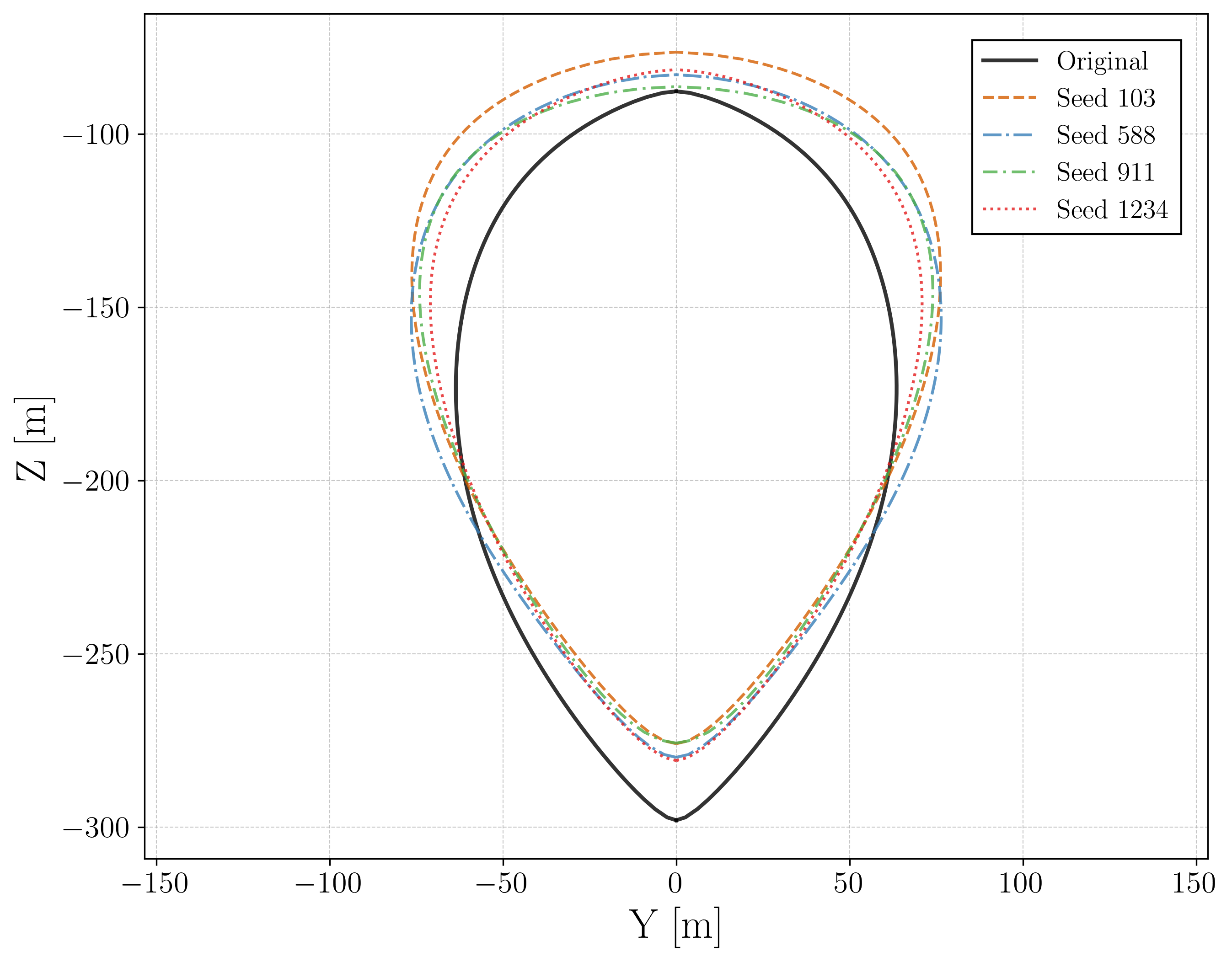} 
        \caption{Extreme Stern (1/20)}
        \label{fig:seed_01}
    \end{subfigure}
    \hfill
    \begin{subfigure}[b]{0.49\textwidth}
        \centering
        \includegraphics[width=\textwidth]{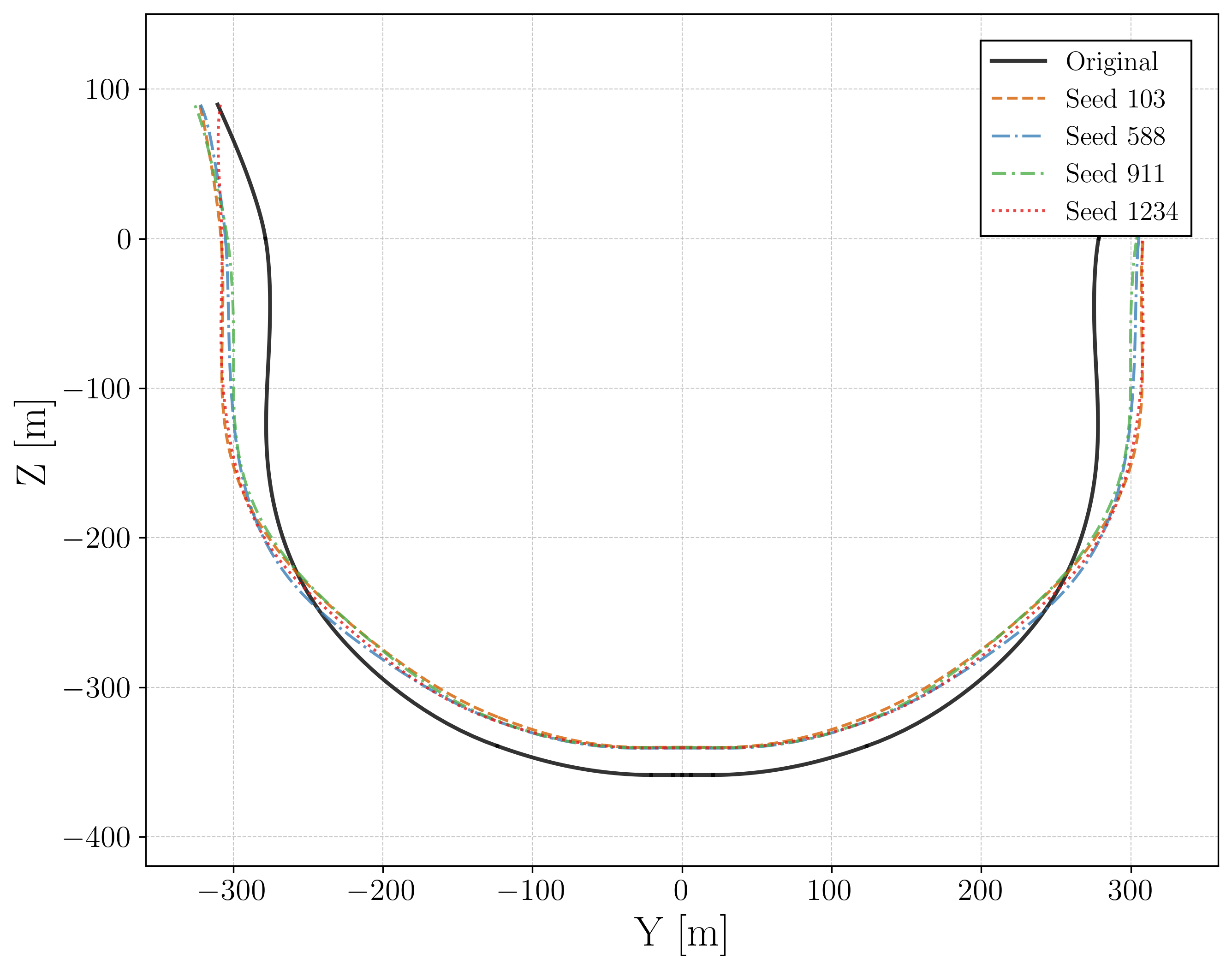} 
        \caption{Aft Section (2/20)}
        \label{fig:seed_02}
    \end{subfigure}

    \vspace{10pt} 

    \begin{subfigure}[b]{0.49\textwidth}
        \centering
        \includegraphics[width=\textwidth]{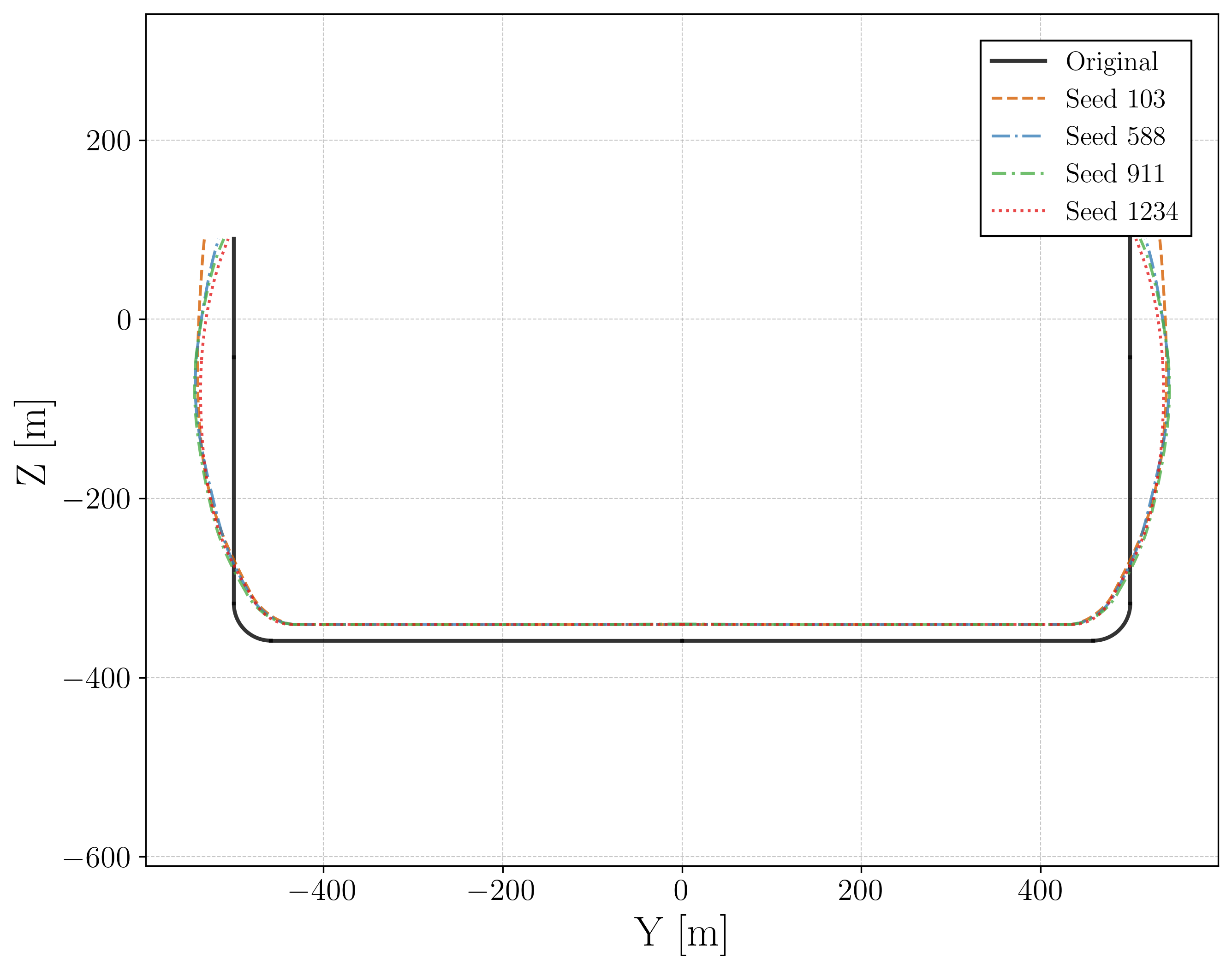} 
        \caption{Midship Section (11/20)}
        \label{fig:seed_11}
    \end{subfigure}
    \hfill
    \begin{subfigure}[b]{0.49\textwidth}
        \centering
        \includegraphics[width=\textwidth]{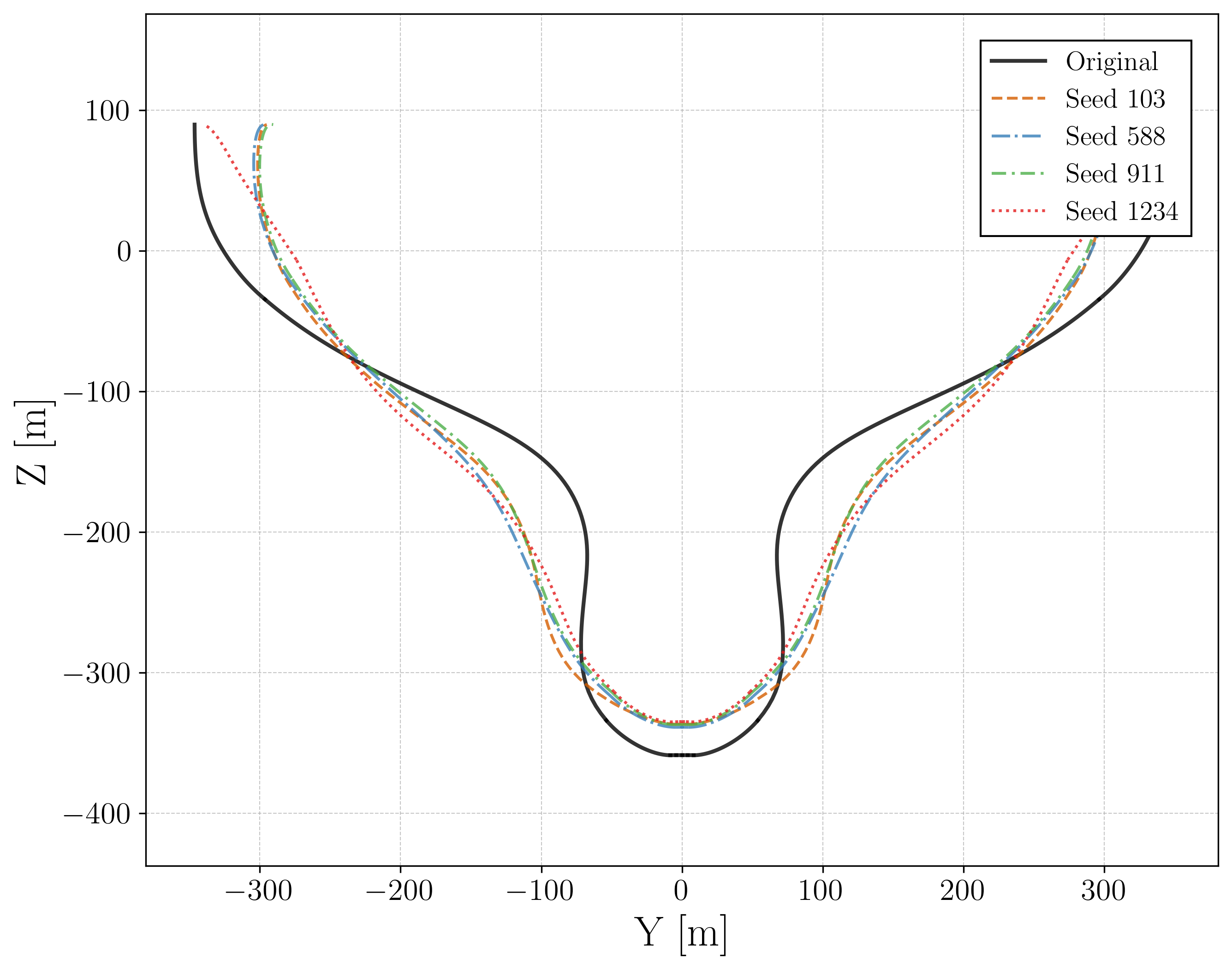} 
        \caption{Bow Section (18/20)}
        \label{fig:seed_18}
    \end{subfigure}
    
    \caption{Robustness analysis of the neural parametrization with respect to stochastic initialization. Comparison of cross-sections at various longitudinal stations for multiple random seeds.}
    \label{fig:seeds_sections}
\end{figure*}

Despite the different starting configurations of the MLP ($\mathcal{N}_{\theta_0}$), the solver consistently converges, exhibiting interesting marginal deviations (\Cref{fig:seeds_sections}). The standard deviation of the final volume error is very low ($\sigma_{vol} \approx 0.08\%$), and the wetted surface area reduction remains consistent across all runs (approx. $2.73\% \pm 0.31\%$). Notably, the draft constraint exhibits the highest stability, with the final minimum vertical coordinate $z_{min}$ varying by less than $4\cdot 10^{-5}$\,\unit{m} across trials. This confirms that the penalty-based barrier function $\mathcal{L}_{depth}$ effectively enforces the hard geometric limit regardless of the initial exploration path. This is due also to the highest value of $\lambda$ in the composite loss function (\Cref{tab:hyperparameters}). The longitudinal centroid shift shows minor fluctuations ($0.0$\,\unit{m} to $2\cdot 10^{-5}$\,\unit{m}), suggesting that while the macroscopic hydrostatic properties are statistically equivalent, the network identifies slightly different local deformation maps to satisfy the global integral constraints, described in \Cref{tab:different_seeds}.

The $N=20$ geometric sections $\Gamma_k$ found in \Cref{fig:seeds_sections} are derived by computing the intersection between the NURBS B-Rep surface $\partial \Omega$ and a set of equidistant transverse planes $\Pi_k$:
\begin{equation}
    \Gamma_k = \partial \Omega \cap \Pi_k, \quad \Pi_k = \{ \mathbf{x} \in \mathbb{R}^3 \mid x = x_k \}
\end{equation}
where $x_k$ is linearly spaced along the longitudinal axis of the hull. 

\subsection{Sensitivity Analysis of Regularization}

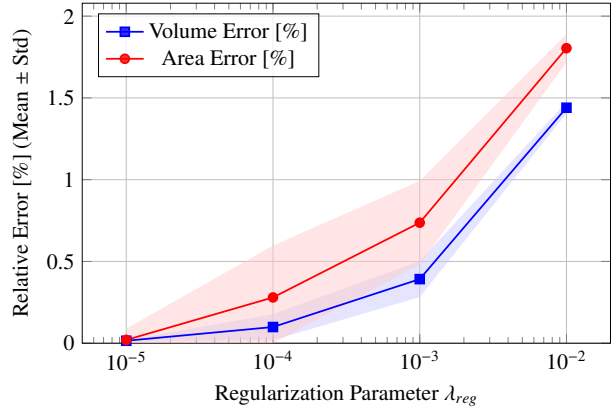
\begin{figure}[htbp]
    \centering
    \resizebox{0.5\textwidth}{!}{
        \begin{tikzpicture}
    \begin{semilogxaxis}[
        width=10cm,
        height=7cm,
        xlabel={Regularization Parameter $\lambda_{reg}$},
        ylabel={Relative Error [\%] (Mean $\pm$ Std)},
        grid=major,
        legend pos=north west,
        legend style={fill=white, fill opacity=0.8, draw opacity=1, text opacity=1},
        cycle list name=color list,
        ymin=0
    ]
    
    \addplot[name path=vol_top, draw=none, forget plot] coordinates {
        (1e-5, 0.018) (1e-4, 0.178) (1e-3, 0.502) (1e-2, 1.475)
    };
    \addplot[name path=vol_bot, draw=none, forget plot] coordinates {
        (1e-5, 0.012) (1e-4, 0.019) (1e-3, 0.282) (1e-2, 1.405)
    };
    \addplot[blue!20, fill opacity=0.5, forget plot] fill between[of=vol_top and vol_bot];
    
    \addplot[color=blue, mark=square*, thick] coordinates {
        (1e-5, 0.015) (1e-4, 0.099) (1e-3, 0.392) (1e-2, 1.440)
    };
    \addlegendentry{Volume Error [\%]}
    
    \addplot[name path=area_top, draw=none, forget plot] coordinates {
        (1e-5, 0.089) (1e-4, 0.597) (1e-3, 0.993) (1e-2, 1.890)
    };
    \addplot[name path=area_bot, draw=none, forget plot] coordinates {
        (1e-5, 0.000) (1e-4, 0.000) (1e-3, 0.481) (1e-2, 1.717)
    };
    \addplot[red!20, fill opacity=0.5, forget plot] fill between[of=area_top and area_bot];
    
    \addplot[color=red, mark=*, thick] coordinates {
        (1e-5, 0.020) (1e-4, 0.280) (1e-3, 0.737) (1e-2, 1.804)
    };
    \addlegendentry{Area Error [\%]}
    
    \end{semilogxaxis}
\end{tikzpicture}
    }    
    \caption{Sensitivity analysis of the geometric optimization performance with respect to the regularization weight $\lambda_{reg}$ for the KVLCC2 hull. Solid lines represent the mean relative error [\%] for volume (blue squares) and wetted surface area (red circles) aggregated over 20 stochastic initializations. Shaded regions denote the $\pm 1$ standard deviation interval, illustrating the convergence of the neural parametrization.}
    \label{fig:lambda_trend}
\end{figure}

The sensitivity of the optimization framework to the regularization weight $\lambda_{reg}$ was investigated by varying the parameter across four orders of magnitude, from $10^{-5}$ to $10^{-2}$, averaging the results among 20 different seeds. 
\Cref{tab:kvlcc2_regularisation} and \Cref{fig:lambda_trend} present the resulting metrics, revealing a fundamental trade-off between the satisfaction of hydrodynamic constraints and the preservation of the initial hull geometry.

\begin{table*}[htbp]
\centering
\caption{Hull Optimization: Regularization sensitivity analysis (Mean $\pm$ SD, $n=20$ for every $\lambda$ value)}
\label{tab:kvlcc2_regularisation}
\footnotesize
\setlength{\tabcolsep}{1.4pt}
\begin{tabular*}{\textwidth}{@{\extracolsep{\fill}} l cc cccc @{}}
\toprule
\textbf{Metric} & \textbf{Initial} & \textbf{Target} & $\lambda=10^{-5}$ & $\lambda=10^{-4}$ & $\lambda=10^{-3}$ & $\lambda=10^{-2}$ \\
\midrule
Volume [\unit{m^3}] & 0.80165 & 0.80165 & $0.80153 \pm 0.00002$ & $0.80086 \pm 0.00064$ & $0.79850 \pm 0.00088$ & $0.79010 \pm 0.00028$ \\
Vol. Error [\%] & -- & 0.0 & $-0.015 \pm 0.003$ & $-0.099 \pm 0.079$ & $-0.392 \pm 0.110$ & $-1.440 \pm 0.035$ \\
\midrule
Area [\unit{m^2}] & 4.12894 & 4.00507 & $4.00588 \pm 0.00275$ & $4.01630 \pm 0.01268$ & $4.03457 \pm 0.01025$ & $4.07731 \pm 0.00347$ \\
Area Err. [\%] & -- & 0.0 & $+0.020 \pm 0.069$ & $+0.280 \pm 0.317$ & $+0.737 \pm 0.256$ & $+1.804 \pm 0.087$ \\
$\Delta$ Area [\%] & -- & -3.0 & $-2.980 \pm 0.067$ & $-2.728 \pm 0.307$ & $-2.285 \pm 0.248$ & $-1.250 \pm 0.084$ \\
\midrule
Max Draft [\unit{m}] & -0.35873 & -0.34079 & $-0.34085 \pm 0.00002$ & $-0.34090 \pm 0.00003$ & $-0.34102 \pm 0.00005$ & $-0.34131 \pm 0.00005$ \\
Draft Dev. [\%] & -- & 0.0 & $-0.018 \pm 0.004$ & $-0.032 \pm 0.010$ & $-0.068 \pm 0.014$ & $-0.151 \pm 0.014$ \\
\midrule
Shift $L_2$ [\unit{m}] & -- & 0.0 & $0.01761 \pm 0.00078$ & $0.01566 \pm 0.00112$ & $0.01334 \pm 0.00061$ & $0.00980 \pm 0.00006$ \\
Shift $x$ [\unit{m}] & -- & 0.0 & $0.00000 \pm 0.00000$ & $0.00001 \pm 0.00001$ & $0.00003 \pm 0.00001$ & $0.00014 \pm 0.00001$ \\
\bottomrule
\end{tabular*}
\end{table*}

Increasing $\lambda_{reg}$ progressively stiffens the deformation field, reducing the mean centroid shift ($L_2$ norm) from $0.018$\,\unit{m} at $\lambda_{reg} = 10^{-5}$ to $0.010$\,\unit{m} at $\lambda_{reg} = 10^{-2}$. However, this geometric rigidity severely compromises the solver's ability to reach the target objectives. At the highest regularization level ($\lambda_{reg} = 10^{-2}$), the mean displaced volume error degrades to $-1.44$\%, a value incompatible with the precision required for hydrodynamic design. In contrast, regimes with $\lambda_{reg} \le 10^{-4}$ maintain volume and wetted area relative errors below $0.3$\%.

\begin{figure}[htbp]
    \centering
    \begin{subfigure}[b]{1\textwidth}
        \centering\includegraphics[width=\textwidth]{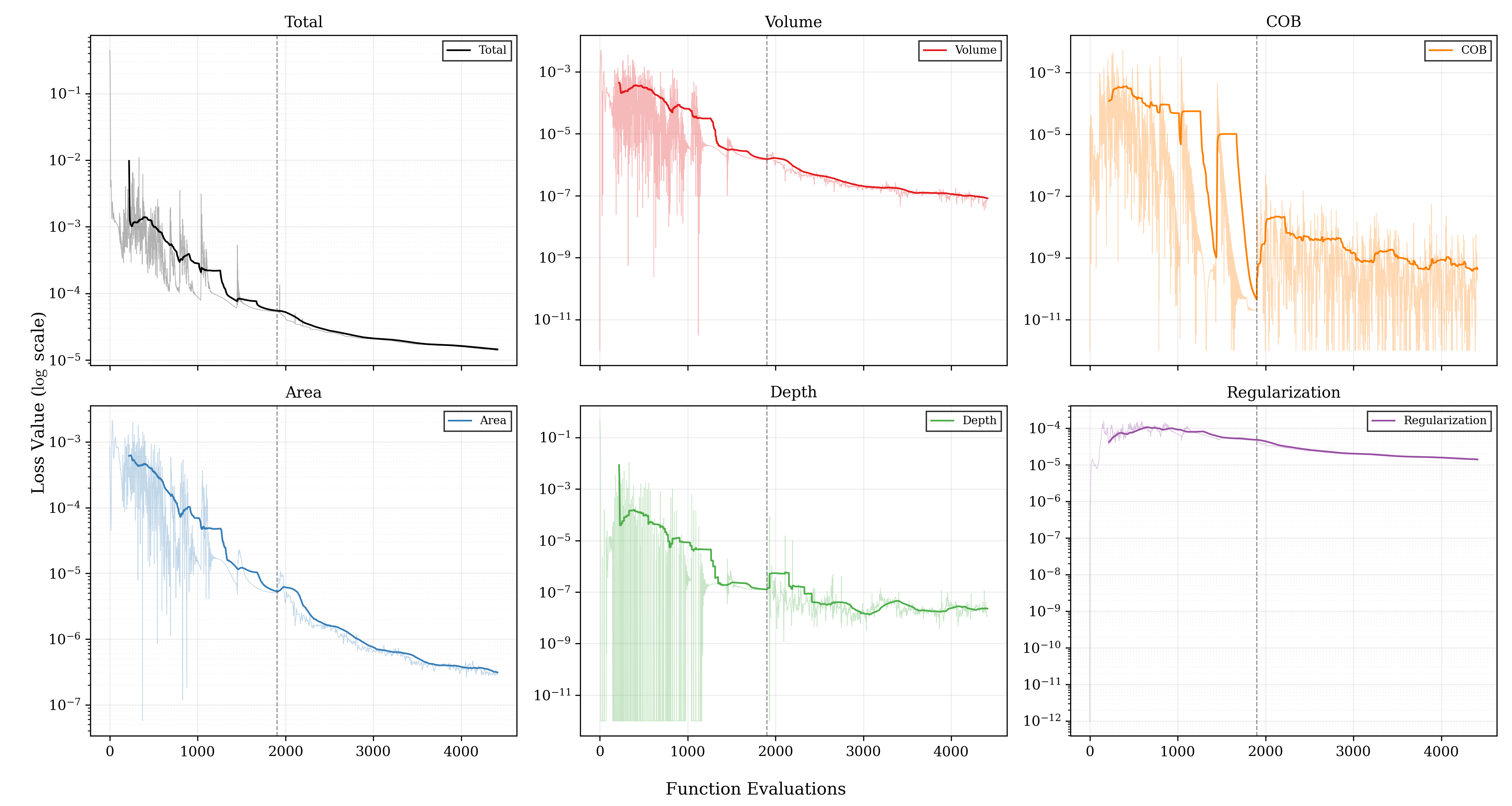} 
        \caption{Example of Low Regularization ($\lambda = 10^{-5}$). Physics driven convergence.}
        \label{fig:loss_low}
    \end{subfigure}
    
    \begin{subfigure}[b]{1\textwidth}
        \centering
        \includegraphics[width=\textwidth]{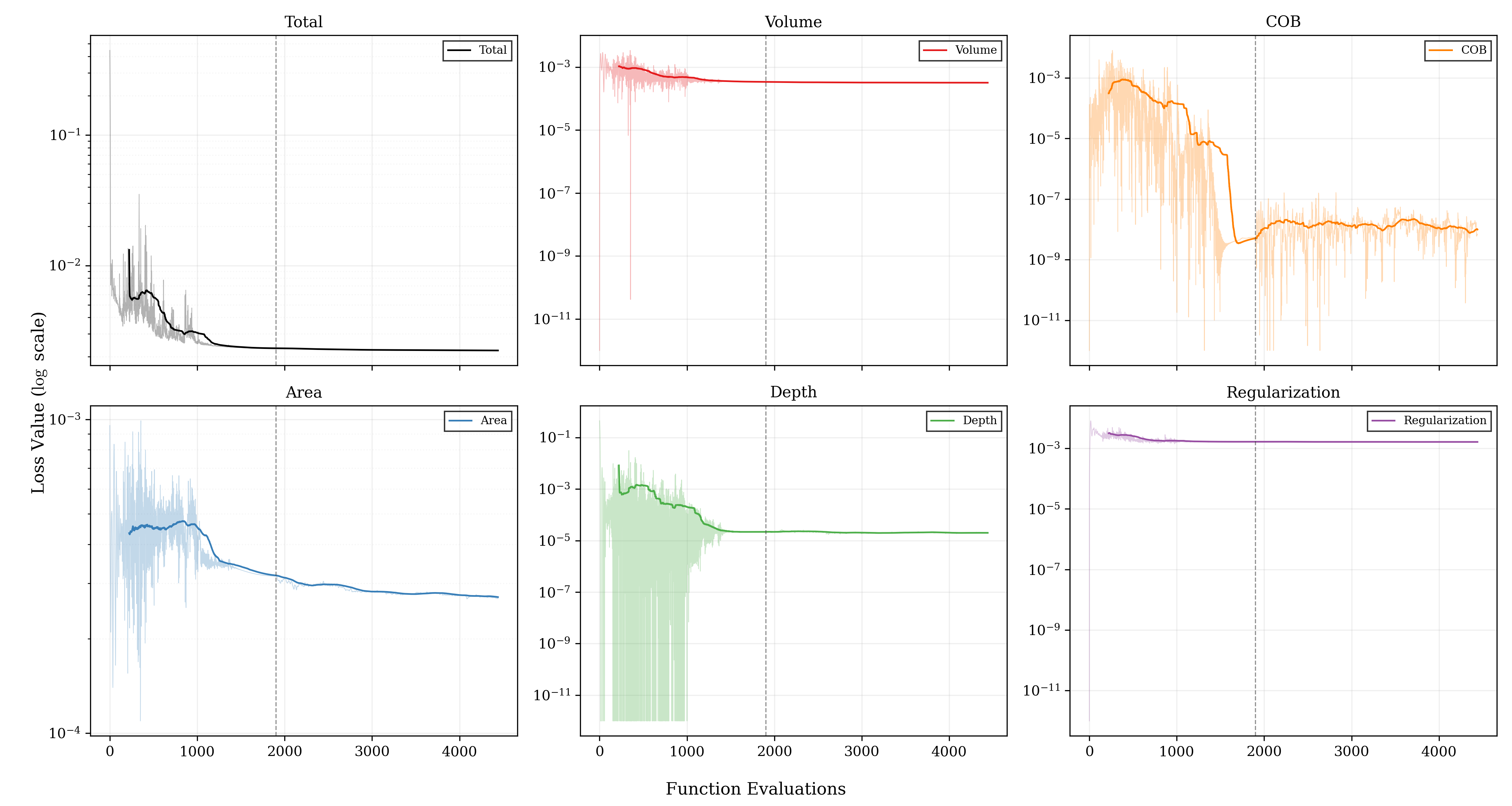} 
        \caption{Example of High Regularization ($\lambda = 10^{-2}$). Regularization dominated regime.}
        \label{fig:loss_high}
    \end{subfigure}
    
    \caption{Comparison of optimisation dynamics for two representative regularization regimes (Seed 911). The vertical dashed line marks the switch from Adam to L-BFGS optimizer. (a) Low $\lambda_{reg}$ allows the second-order optimizer to effectively minimize physical residuals. (b) High $\lambda_{reg}$ leads to stagnation as the optimization landscape is dominated by the regularization penalty.}
    \label{fig:loss_comparison}
\end{figure}

The underlying optimization dynamics are visualized in  \Cref{fig:loss_comparison}, which contrasts the optimisation histories of the limiting regimes. In a low regularization case ($\lambda_{reg} = 10^{-5}$), the transition from the  AdamW \cite{adamw} optimizer to the L-BFGS \cite{L-BFGS} algorithm triggers a reduction in total residuals of approximately two orders of magnitude, indicating successful convergence to the physical constraints. Conversely, the stiff regime ($\lambda_{reg} = 10^{-2}$) is dominated by the regularization gradient throughout the process. Based on these findings, $\lambda_{reg} = 10^{-4}$ is identified a possible optimal compromise to ensure negligible physical errors while preventing mesh degeneracy.

\begin{figure*}[htbp]
    \centering
    \begin{subfigure}[b]{0.49\textwidth}
        \centering
        \includegraphics[width=\textwidth]{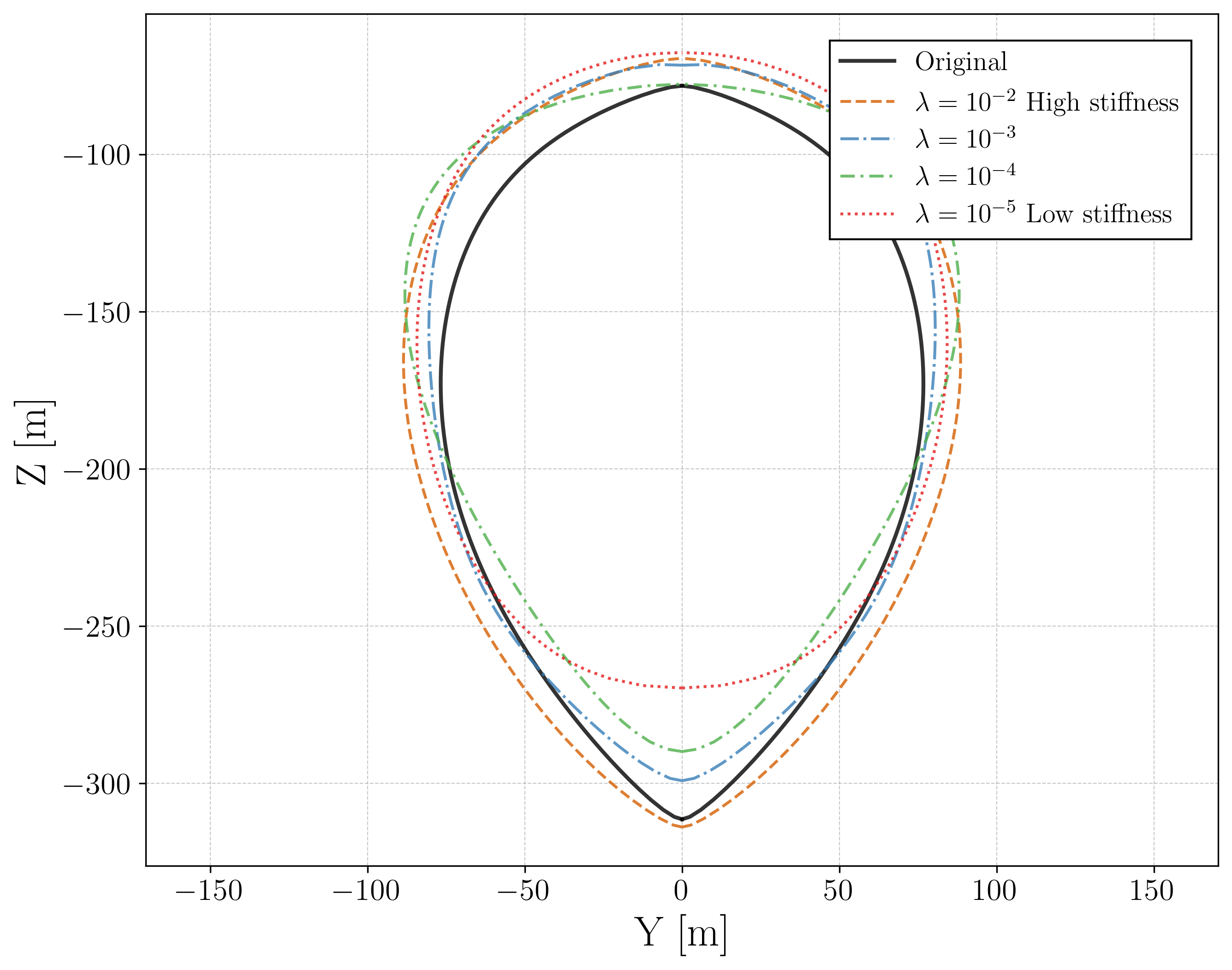} 
        \caption{Extreme Stern (1/20)}
        \label{fig:sec_01}
    \end{subfigure}
    \hfill
    \begin{subfigure}[b]{0.49\textwidth}
        \centering
        \includegraphics[width=\textwidth]{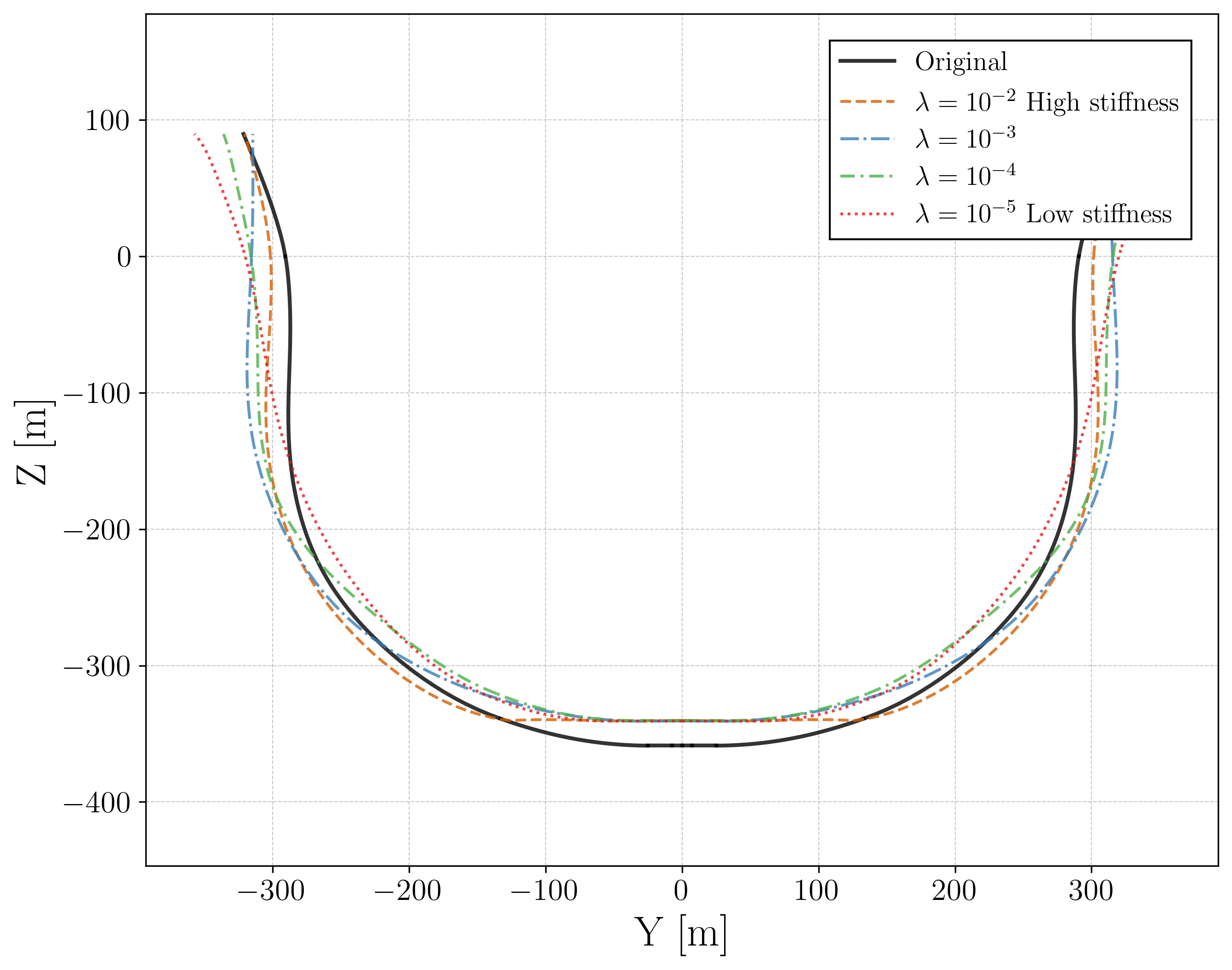} 
        \caption{Aft Section (2/20)}
        \label{fig:sec_02}
    \end{subfigure}

    \vspace{10pt} 

    \begin{subfigure}[b]{0.49\textwidth}
        \centering
       \includegraphics[width=\textwidth]{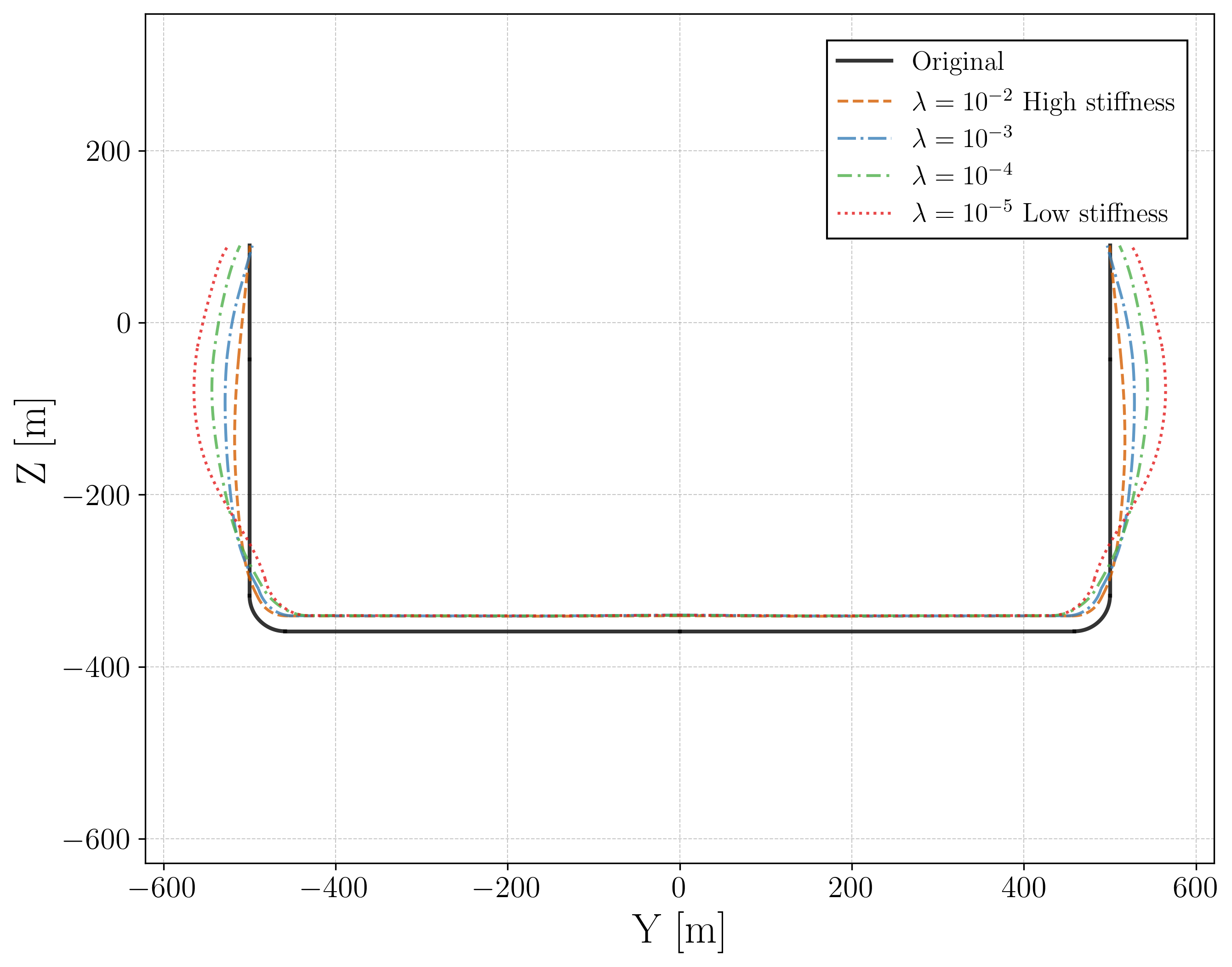} 
        \caption{Midship Section (11/20)}
        \label{fig:sec_mid}
    \end{subfigure}
    \hfill
    \begin{subfigure}[b]{0.49\textwidth}
        \centering
        \includegraphics[width=\textwidth]{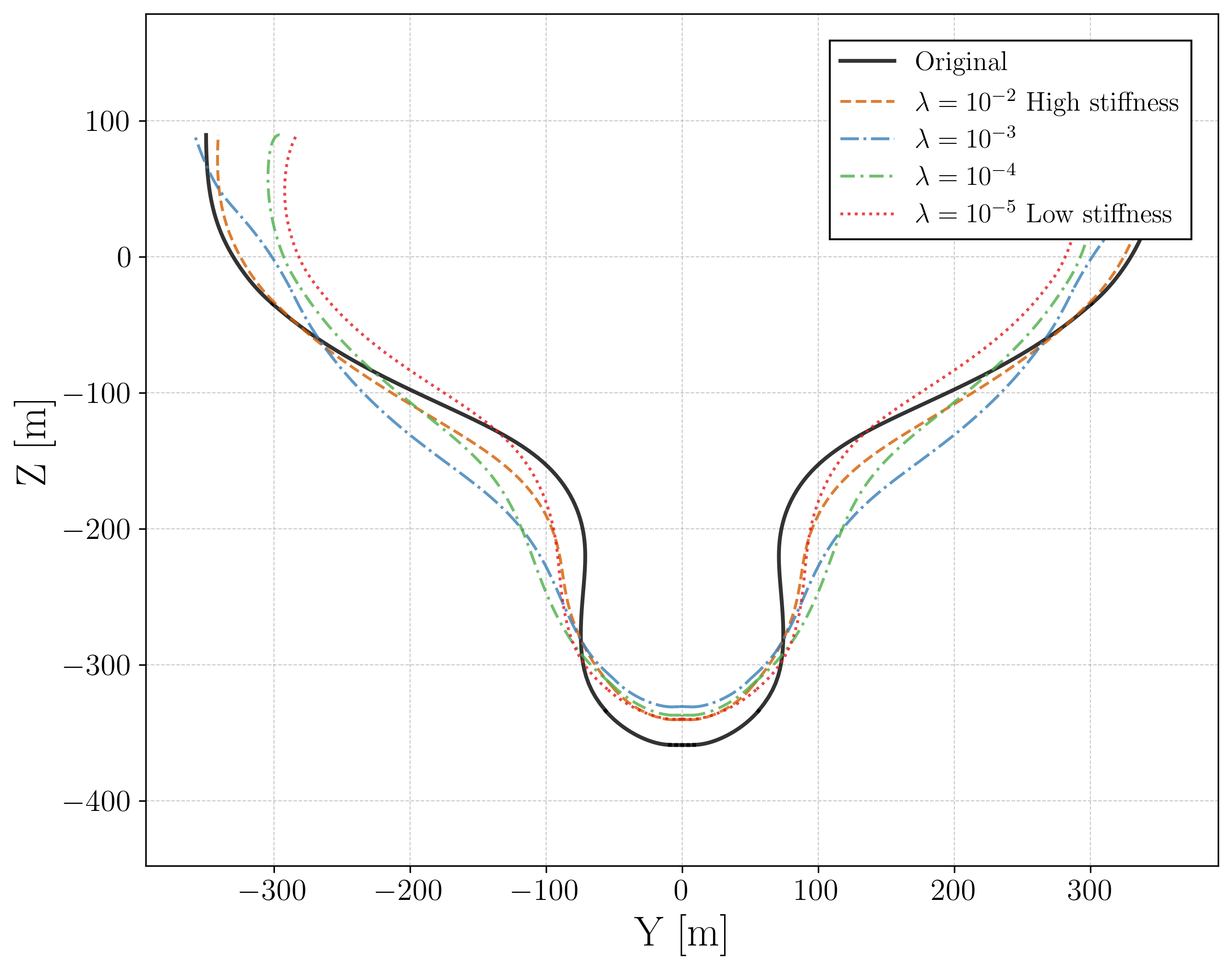} 
        \caption{Bow Section (18/20)}
        \label{fig:sec_bow}
    \end{subfigure}
    
    \caption{\textit{Sensitivity of geometric deformation to regularization stiffness.} 
    Cross-sections at four representative stations. 
    \textbf{(c)} At midship, the \textit{Low Stiffness} regime (orange, $\lambda_{reg}=10^{-5}$) enables the necessary beam expansion to compensate for the draft reduction ($Z_{min}$ constraint), ensuring volume conservation within numerical tolerance. 
    In contrast, the \textit{High Stiffness} regime (red, $\lambda_{reg}=10^{-2}$) locks the geometry to the original manifold, leading to the volume violations. The deformation field is global: mass redistribution is evident along the entire hull to maintain the longitudinal center of buoyancy ($x_{COB}$). The random seed used is 911.}
    \label{fig:sensitivity_sections}
\end{figure*}

The numerical trade-off observed is physically shown by the cross-sectional analysis in \Cref{fig:sensitivity_sections}. At the midship section (\Cref{fig:sec_mid}), satisfying the draft constraint ($z_{min}$) forces a vertical displacement of the keel. To conserve the total volume, the solver must redistribute the displaced mass by expanding the hull beam.
As observed in the cross-sectional analysis (\Cref{fig:sec_mid}), the unconstrained volume compensation alters the verticality of the ship sides. While hydrostatically viable for an oil carrier, this geometric deformation could introduce inefficiencies in vessel classes requiring vertical side walls, such as container ships. 
Further dimensional constraints, analogous to the enforced draft barrier, can be directly integrated into the optimization objective to preserve specific structural features based on the designer's operational requirements.

The high regularization regime ($\lambda_{reg}=10^{-2}$) imposes excessive stiffness on the control points, preventing this lateral expansion. This results in the significant volume loss discussed above. Conversely, the low regularization regime ($\lambda_{reg}=10^{-5}$) grants the neural network sufficient flexibility to identify a non-rigid compensation mechanism, widening the hull sections while maintaining surface fairness. Furthermore, Figures \ref{fig:sec_01}, \ref{fig:sec_02} and \ref{fig:sec_bow} demonstrate that this deformation is not localized but globally distributed towards the stern and bow, ensuring the preservation of the longitudinal center of buoyancy ($x_{COB}$).


\section{Conclusions}
\label{sec:conclusions}
This work established an end-to-end differentiable framework for the shape optimization of parametric surfaces, implemented within the PINA ecosystem. By leveraging analytical B-Spline derivatives for the computation of the Jacobian determinant, the proposed method enables fast steps for integral constraints problems. This optimisation formulation avoids the geometrical discretization overhead and gradient inconsistencies inherent to traditional mesh-based analysis workflows, and more flexibility than traditional Free Form Deformations standard approaches.

A key contribution is the introduction of a neural deformation strategy, where a Multi-Layer Perceptron (MLP) parametrizes the global displacement field of the control net. Comparative benchmarks against direct control point optimization demonstrate that the MLP architecture functions as an implicit regularizer. This regularization is critical for preserving $C^0$ continuity across multi-patch interfaces, effectively preventing the geometric degradation observed in direct approaches.

Numerical validation on the hull geometry confirmed the framework's efficacy in handling competing geometric constraints. The solver successfully enforced a 5\% draft reduction constraint and achieved an approximate 2.8\% average reduction in wetted surface area (against a 3\% target), while maintaining the longitudinal geometry centroid. The method achieved volume conservation with a mean relative error below 0.1\%, producing geometries with interface continuity. 
Sensitivity analysis identified a regularization regime ($\lambda_{reg} \approx 10^{-4}$)  ensuring convergence to valid design candidates.

The intrinsic generality of the proposed framework allows for its extension to heterogeneous engineering domains, provided that the integral geometric constraints and boundary conditions are appropriately adapted to the specific physical problem.

Future research can focus on two primary directions. First, the framework can be extended to support trimmed NURBS surfaces to address arbitrary industrial topologies defined by non-tensor product boundaries. Second, the pipeline can incorporate PINN solvers to resolve Partial Differential Equations (PDEs) directly on the evolving manifold. Leveraging the PINA infrastructure, this extension aims to transition from geometric constraint satisfaction to functional optimization, targeting physical performance metrics such as hydrodynamic drag or structural stress. In addition, the actual discrete integration can be replaced by the analytical formulation of B-Spline polynomial integration.


\section*{Acknowledgements}

F.T., G.A.D., N.D. and G.R. acknowledge the support provided by the European Union -- NextGenerationEU, in the framework of the iNEST -- Interconnected Nord-Est Innovation Ecosystem (iNEST ECS00000043 – CUP G93C22000610007) project and its CC5 Young Researchers initiative. G.C. acknowledges the support provided by the European Union -- NextGenerationEU, within the framework of the National Recovery and Resilience Plan (PNRR), Mission 4, Component 2, Investment 3.3, through a PhD scholarship at Scuola Internazionale Superiore di Studi Avanzati di Trieste (CUP G93C24000530001), developed in collaboration with Fincantieri S.p.A. G.R. acknowledges the support provided by PRIN "FaReX - Full and Reduced order modelling of coupled systems: focus on non-matching methods and automatic learning" project.
The views and opinions expressed are solely those of the authors and do not necessarily reflect those of the European Union, nor can the European Union be held responsible for them. G.A.D. and G.R. would like to acknowledge INdAM–GNCS.


\bibliographystyle{elsarticle-num} 
\bibliography{references}

\end{document}